\documentclass{amsart}

\theoremstyle{remark}

\numberwithin{equation}{section}

\usepackage{graphicx} 
\usepackage{extarrows}
\usepackage{latexsym}
\usepackage{amsmath}
\usepackage{amssymb}
\usepackage{amsfonts}
\usepackage{verbatim}
\usepackage{mathrsfs}
\usepackage[modulo]{lineno} 
\usepackage{color}
\usepackage{xcolor}
\usepackage[colorlinks,citecolor=blue,urlcolor=blue]{hyperref}
\usepackage{soul}

\newtheorem{theorem}{Theorem}[section]
\newtheorem{remark}{Remark}[section]
\newtheorem{assumption}{Assumption}[section] 
\newtheorem{proposition}{Proposition}[section]
\newtheorem{lemma}{Lemma}[section] 

\newcommand{\tr}{\top}

\newcommand{\ee}{\mathbb E}

\newcommand{\pp}{\mathbb P}
\newcommand{\nn}{\mathbb N}
\newcommand{\rr}{\mathbb R}

\newcommand{\BB}{\mathcal B}
\newcommand{\CC}{\mathcal C}
\newcommand{\LL}{\mathcal L}

\newcommand{\PP}{\mathcal P}

\newcommand{\OOO}{\mathscr O}

\newcommand{\<}{\langle}
\renewcommand{\>}{\rangle}
\allowdisplaybreaks \allowdisplaybreaks[4]

\begin{document}

\title[Ergodicity and Estimates for Tamed Schemes of superlinear SPDEs]
{Numerical Ergodicity and Optimal Strong Error Estimates for a Class of Novel Tamed Schemes to Superlinear SPDEs}
 
\author{Zhihui LIU}
\address{Department of Mathematics \& National Center for Applied Mathematics Shenzhen (NCAMS) \& Shenzhen International Center for Mathematics, Southern University of Science and Technology, Shenzhen, 518055, P.R.China }
\email{liuzh3@sustech.edu.cn (Corresponding author)}
 
\author{Jie SHEN}
\address{Distinguished Professor Emeritus, Purdue University; Chair Professor and Dean, School of Mathematical Science, Eastern Institute of Technology, Ningbo, 315200, P.R.China}
\email{jshen@eitech.edu.cn}

\thanks{}


 

\begin{abstract} 
We construct a class of novel tamed schemes for superlinear stochastic partial differential equations (SPDEs), including the stochastic Allen--Cahn equation driven by either multiplicative or additive noise. The schemes preserve the same Lyapunov structure as the original system, and we rigorously establish their longtime unconditional stability. Furthermore, we prove that the corresponding Galerkin-based fully discrete tamed schemes inherit the unique ergodicity of the underlying SPDEs and achieve optimal strong convergence rates in both the multiplicative and additive noise cases.
\end{abstract}

\maketitle

\section{Introduction}

The longtime numerical behavior --- in particular, the numerical ergodicity --- of superlinear SPDEs whose coefficients violate the classical Lipschitz continuity and linear growth assumptions remains an intriguing and challenging problem. As a fundamental longtime property, ergodicity characterizes the equivalence between temporal and spatial averages (the so-called ergodic limit) for Markov processes and chains generated by SPDEs and their numerical discretizations, respectively. This concept is pivotal in diverse fields such as quantum mechanics, fluid dynamics, and financial mathematics \cite{DZ96}, \cite{HW19}.

In practical applications, a central goal is to compute the ergodic limit, i.e., the mean of a given function with respect to the equilibrium distribution. This distribution is also known as the invariant measure. However, obtaining an explicit expression for the invariant measure of an infinite-dimensional stochastic system is typically infeasible. Consequently, it is generally impossible to compute the ergodic limit directly for nonlinear SPDEs driven by multiplicative noise. This challenge has motivated a substantial body of research over the past decade, leading to the development of numerical algorithms that preserve the ergodicity of the original system and efficiently approximate its ergodic limit.

In this paper, we investigate the superlinear SPDE
\begin{align}\label{see-fg} 
& {\rm d} X(t, \xi)  
=(\Delta X(t, \xi)+f(X(t, \xi))) {\rm d}t
+g(X(t, \xi)) {\rm d}W(t, \xi),  
\end{align}
$(t, \xi) \in \rr_+\times \OOO$, with homogeneous Dirichlet boundary condition $X(t, \xi)=0$, $(t, \xi) \in \rr_+\times \partial \OOO$, and initial condition $X(0, \cdot)=X_0(\cdot)$. Here, $\OOO \subset \rr^d$ ($d=1,2,3$) is a bounded open set with a piecewise smooth boundary, $f,g: \rr \to \rr$ are assumed to be monotone-type with certain polynomial growth, and $W$ is an infinite-dimensional Wiener process (cf. Section \ref{sec2} for more details).
 
Eq. \eqref{see-fg} includes the following stochastic Allen--Cahn equation (SACE), arising from phase transition in materials science by stochastic perturbation, as a special case:
\begin{align} \label{ac}
{\rm d} X(t, \xi)= \Delta X(t, \xi)  {\rm d}t + \frac1{\epsilon^2} (X(t, \xi)-X(t, \xi)^3) {\rm d}t + g(X(t, \xi)) {\rm d}W(t, \xi), 
\end{align}
 where the positive index $\epsilon \ll 1$ is the interface thickness; see, e.g., \cite{BGJK23}, \cite{BP24}, \cite{CH19}, \cite{FLZ17}, \cite{LQ20}, \cite{LQ21}, \cite{OPW23} and references therein.  
Our aim is to construct efficient numerical integrators for Eq. \eqref{see-fg}, whose solution is uniquely ergodic with respect to an equilibrium distribution, that satisfy the following two properties: 
\begin{enumerate}
\item[(1)]
exhibit ergodicity with respect to approximate equilibrium distributions of Eq. \eqref{see-fg} on infinite time intervals; 

\item[(2)]
converge strongly with the optimal rate to the solutions of Eq. \eqref{see-fg} on any finite time intervals. 
\end{enumerate}   

Significant progress has been made in the design and analysis of numerical approximations for ergodic limits in finite-dimensional settings, particularly for stochastic ordinary differential equations \cite{LV22}, \cite{LL24}, \cite{LWWZ25}, \cite{LW26a}, \cite{LW26}, \cite{MSH02}. In contrast, the construction and analysis of numerical ergodic limits for SPDEs, even those driven by additive noise, are still in their early stages. 

The backward Euler method and its Galerkin-based full-discretizations satisfy properties (1) and (2) above; however, they often entail high computational costs. For instance, \cite{Liu26} discusses the exponential ergodicity of SACE driven by multiplicative trace-class noise with a moderately thick interface. In contrast, \cite{Liu24b} establishes the unique ergodicity in the multiplicative white noise case, regardless of the interface thickness. Additionally, \cite{Bre22}, \cite{BK17}, \cite{BV16}, \cite{CHS21}, \cite{WC26} explore approaches for approximating the invariant measures of SPDEs driven by additive white or colored noise. 

Schemes that treat the nonlinearity explicitly are often preferred to reduce computational costs. 
To construct such explicit schemes that can inherit the unique ergodicity of Eq. \eqref{see-fg}, as in the tamed methods studied for the finite-dimensional case (e.g.,  \cite{HJK12}, \cite{LNSZ24}, \cite{NNZ25}), the numerical Lyapunov structure plays a crucial role. However, \cite{BHJKLS19} demonstrated that the classical Euler--Maruyama (EM) scheme and its Galerkin-based full-discretizations, when applied to Eq. \eqref{see-fg} with superlinear growth coefficients, lead to blow-up in the $p$-th moment for all $p \geq 2$. In particular, the square function (in certain norm), although a natural Lyapunov functional for the monotone equation \eqref{see-fg} (which can be verified by applying It\^o formula and a coupled coercivity condition, in combination with a compact embedding argument in infinite-dimensional setting; see Assumption \ref{ap} and \cite{Liu26} for details), is not a suitable Lyapunov functional for the EM scheme.

To address this issue, we propose a family of novel explicit-in-nonlinearity tamed schemes that achieve the two main objectives (1) and (2) for superlinear SPDEs \eqref{see-fg} driven by the more subtle and challenging multiplicative noise. Compared to existing tamed schemes in the literature, such as those in \cite{Bre22}, \cite{HS23}, \cite{WC26}, our schemes are more efficient to implement because they tame the nodal values rather than their norms at each step. Moreover, our carefully chosen taming functions (in \eqref{bs-tau}) preserve the coupled monotonicity (in \eqref{coe-tau}), thereby maintaining the Lyapunov structure that conventional taming disrupts. 

A key step of our approach is the rigorous analysis of the longtime unconditional stability of the temporal-tamed scheme and its Galerkin-based full discretizations. This analysis demonstrates that these fully discrete tamed schemes preserve the same Lyapunov structure as the original superlinear equation \eqref{see-fg}. Combined with the equivalence of transition probabilities under a certain nondegeneracy condition, these tamed schemes are shown to inherit the unique ergodicity of \eqref{see-fg}.

On the other hand, as pointed out in \cite{LW26a}, even for Langevin systems, pure sampling methods can achieve (1) but typically do not approximate the exact solution. In contrast, integrators satisfy (2), yet they often diverge on infinite time intervals or become ergodic with respect to a different equilibrium distribution. Therefore, it is crucial to assess the impact of these tamed methods on the dynamics of \eqref{see-fg} through strong error estimates. 

To address the second objective and to derive higher-order strong convergence rates in the additive noise case (see Theorem \ref{tm-err-}), we point out that, unlike previous strong error analyses in \cite{BGJK23,FLZ17,Liu26,LQ20,LQ21,QW19}, our proof requires the \(L^s\)-regularity (with \(s>2\)) of the exact solution, along with the \(\dot H^1\)-regularity of the auxiliary process associated with the tamed scheme. Furthermore, in the 3D case, such regularity cannot be controlled by the \(\dot {H} ^1 \)-norm alone (see Remark~\ref{rk-reg}), whereas it is not required in the multiplicative noise case (see Theorem \ref{tm-err}).

Our main contributions are as follows:
\begin{itemize}
\item the construction of a class of computationally efficient tamed schemes for superlinear SPDEs that preserve the Lyapunov stucture and unique ergodicity;
\item a rigorous proof that the solutions of these new schemes are unconditionally stable over any time interval;
\item the establishment of strong convergence of these tamed schemes with optimal spatial and temporal orders for both the multiplicative and additive noise cases.
\end{itemize}
To the best of our knowledge, this is the first construction of explicit-in-nonlinearity schemes that preserve the unique ergodicity of superlinear SPDEs driven by either additive or multiplicative noise. We emphasize that our results hold in one, two, and three space dimensions, thereby significantly improving the finite-time stability results in prior work \cite{HS23} for another type of tamed scheme and removing the dimension restrictions imposed therein, while requiring weaker regularity assumptions on the exact solution of Eq.~\eqref{see-fg}. We also note that, since the first version of this paper appeared, our taming idea has been applied in \cite{JW25} to obtain longtime weak error estimates for SACE driven by additive noise and in \cite{CL26} and \cite{CLW26} to establish uniform-in-time strong and weak/ergodic error estimates for Eq. \eqref{see} with certain additional dissipativity conditions, respectively.

The remainder of the paper is organized as follows. Section~\ref{sec2} introduces the necessary preliminaries and the proposed tamed schemes. In Section~\ref{sec3}, we develop the Lyapunov structure and probabilistic regularity of the tamed schemes with Galerkin-based full-discretizations, from which we establish their unique ergodicity. Sections~\ref{sec4} and~\ref{sec5} present the strong error estimates between these tamed full-discretizations and Eq.~\eqref{ac} in the multiplicative and additive noise cases, respectively.

\section{Preliminaries}
\label{sec2}

In this section, we introduce some notations and main assumptions, then propose the tamed schemes. 
  
  \subsection{Notations}
  \label{sec2.1}

For two constants $a$ and $b$, we denote $a \wedge b := \min\{a, b\}$.   
For $p \in [1, \infty]$, we denote by $(L^p(\OOO), \|\cdot\|_{L^p_\xi})$, $(L^p(0,T), \|\cdot\|_{L^p_t})$, and $(L^p(\Omega), \|\cdot\|_{L^p_\omega})$ the usual real-valued Lebesgue spaces on $\OOO$, $(0,T)$ (for a fixed $T \in (0,\infty)$), and $\Omega$, respectively. When necessary, we also use mixed norms such as $\|\cdot\|_{L^p_\omega L^r_t L^l_\xi}$ for $p, r, l \in [1, \infty]$. 

Let $H := \{u \in L^2(\OOO) : u|_{\partial\OOO} = 0\}$ with inner product $\langle \cdot, \cdot \rangle$ and norm $\|\cdot\|$ (or $\|\cdot\|_{L^2_\xi}$). Denote by $A: \operatorname{Dom}(A) \subset H \to H$ the Dirichlet Laplacian on $H$. Then $A$ is the infinitesimal generator of an analytic $C_0$-semigroup $S(\cdot) = e^{A \cdot}$ on $H$, and one can define the fractional powers $(-A)^\theta$ for $\theta \in \mathbb{R}$ of the self-adjoint positive-definite operator $-A$. Let $\dot H^\theta$ be the domain of $(-A)^{\theta/2}$ equipped with the norm $\|\cdot\|_\theta$ given by $\|u\|_\theta := \|(-A)^{\theta/2} u\|$ for $u \in \dot H^\theta$. In particular, $\dot H^0 = H$ and $\dot H^1 = H^1_0 := \{u \in L^2(\OOO) : \nabla u \in L^2(\OOO),\, u|_{\partial\OOO}=0\}$, whose inner product is $\<\cdot, \cdot\>_1:=\<\cdot, \cdot\>+\<\nabla \cdot, \nabla \cdot\>$. The dual space of $V := \dot H^1$ is denoted by $V^* := \dot H^{-1}$, with duality pairing ${}_{-1}\langle \cdot, \cdot \rangle_1$ (with respect to $\langle \cdot, \cdot \rangle$). We use $\operatorname{Id}$ to denote the identity operator on any Hilbert space when no confusion arises. 

Let $U$ be a separable Hilbert space and $Q$ a self-adjoint positive-definite operator on $U$. Set $U_0:= Q^{1/2}U$ and denote by $(\mathcal{L}_2^\theta := HS(U_0; \dot H^\theta), \|\cdot\|_{\mathcal{L}_2^\theta})$ the space of Hilbert--Schmidt operators from $U_0$ to $\dot H^\theta$ for $\theta \in \mathbb{R}_+$. Let $W := \{W(t) : t \ge 0\}$ be a $U$-valued $Q$-Wiener process on the stochastic basis $(\Omega,\mathcal F, \{\mathcal F_t\}_{t \ge 0}, \pp)$, i.e., there exist an orthonormal basis $\{g_m\}_{m=1}^\infty$ of $U$ consisting of eigenvectors of $Q$ with corresponding eigenvalues $\{q_m\}_{m=1}^\infty$, and a sequence of mutually independent 1D Brownian motions $\{\beta_m\}_{m=1}^\infty$ such that $W(t) = \sum_{m=1}^\infty \sqrt{q_m}\, g_m \beta_m(t)$ for $t \ge 0$. We mainly focus on the trace-class case, i.e., $\sum_{m=1}^\infty q_m < \infty$; we refer to \cite{Liu24b} for the numerical ergodicity of an implicit scheme for Eq.~\eqref{see-fg} driven by white noise.

\subsection{Assumptions}

To study the numerical Lyapunov structure, we impose the following coupled coercivity and polynomial growth conditions on the measurable functions $f,g: \mathbb{R} \to \mathbb{R}$ in Eq. \eqref{see-fg}.

\begin{assumption} \label{ap}   
There exist positive constants $L_i$, $i=1,2,3,4$, and $q>0$ such that for any $\xi \in \rr$ and any function $u: \rr \to \rr$, 
\begin{align} 
u(\xi) f(u(\xi))+ |g(u(\xi))|^2 \sum_{m \in \nn_+} q_m |g_m(\xi)|^2 
& \le L_1 - L_2 |u(\xi)|^{q+2}, \label{coe} \\
|f(\xi)| & \le L_3+ L_4 |\xi|^{q+1}. \label{f-grow}
\end{align}  
 \end{assumption}

Throughout, we assume that $q>0$ when $d=1,2$ and $q \in (0,2]$ when $d=3$, so that the following Sobolev embeddings hold (see \cite{LQ21} for the justification of this range of $q$):
\begin{align} \label{emb}
\dot H^1\hookrightarrow L_\xi^{2(q+1)} 
\hookrightarrow L_\xi^2 \hookrightarrow L_\xi^{[2(q+1)]'} \hookrightarrow \dot H^{-1},
\quad d=1,2,3, 
\end{align}
where $[2(q+1)]':=2(q+1)/(2q+1)$ denotes the conjugate exponent of $2(q+1)$.
When $d = 3$, the paper's analysis is restricted to the SACE. 
Then we can define the Nemytskii operators $F: \dot H^1 \rightarrow \dot H^{-1}$ and $G: H \rightarrow \LL_2^0$ associated with $f$ and $g$, respectively, by
\begin{align} 
F(u)(\xi):= f(u(\xi)), & \quad u \in \dot H^1, ~ \xi \in \OOO, \label{df-F}\\
G(u) g_m(\xi):= g(u(\xi)) g_m(\xi), & \quad u \in H,~ m \in \nn_+,~ \xi \in \OOO.
\label{df-G}
\end{align}
   
With these preliminaries, Eq. \eqref{see-fg} is equivalent to the infinite-dimensional stochastic evolution equation
\begin{align} \label{see}
{\rm d} X(t)=(A X(t)+F(X(t))) {\rm d}t+G(X(t)) {\rm d}W, \quad t \ge 0;
\quad X(0)=X_0.
\end{align} 

 \begin{remark} 
Conditions \eqref{coe}-\eqref{f-grow} imply a polynomial growth condition on $g$, i.e., there exist positive constants $L_5$ and $L_6$ such that 
 \begin{align*}
|g(u(\xi))|^2 \sum_{m \in \nn_+} q_m |g_m(\xi)|^2 
\le L_5+L_6 |u(\xi)|^{q+2},
\quad \xi \in \rr, ~  u: \rr \to \rr,
 \end{align*} 
which yields 
$\|G(u)\|_{\LL_2^0} ^2 \le L_5 |\OOO| +L_6 \|u\|_{L_\xi^{q+2}}^{q+2}$, $u \in L_\xi^{q+2}$.
Condition \eqref{coe} also implies the following infinite-dimensional version of the coercive condition: 
\begin{align*} 
_{L_\xi^{q+2}}\<u, F(u)\>_{L_\xi^{(q+2)'}} + \|G(u)\|_{\LL_2^0} ^2 
\le L_1 |\OOO| - L_2 \|u\|_{L_\xi^{q+2}}^{q+2},
\quad u \in L_\xi^{q+2},
 \end{align*} 
 where $|\OOO|$ denotes the volume of $\OOO$ and $ _{L_\xi^{q+2}}\<\cdot, \cdot\>_{L_\xi^{(q+2)'}}$ is the duality pairing between $L_\xi^{q+2}$ and its conjugate space $L_\xi^{(q+2)'}$.
 In the case $\sum_{m \in \nn_+} q_m \|g_m\|_{L_\xi^\infty}^2 < \infty$, 
which holds when $\sup_{m \in \nn_+} \|g_m\|_{L_\xi^\infty} < \infty$, condition \eqref{coe} is equivalent to  
 \begin{align*} 
\xi f(\xi)+ \sum_{m \in \nn_+} q_m \|g_m\|_{L_\xi^\infty}^2 |g(\xi)|^2 
& \le L_1 - L_2 |\xi|^{q+2}, \quad \xi \in \rr.
\end{align*} 
 \end{remark} 

\begin{remark} \label{rk-ap-pol} 
Assumption~\ref{ap} includes certain polynomial drifts (with a negative leading coefficient) and diffusion functions. Indeed, let $k \in \mathbb{N}$ and assume
 \begin{align}\label{bs-ex}
f(\xi)=\sum_{i=0}^{2k+1} a_i \xi^i, \quad
g(\xi)=\sum_{j=0}^{k+1} c_j \xi^j,
\quad \xi \in \rr,
 \end{align} 
 where $a_i, c_j \in \rr$, $i=0,1,\cdots,2k+1$, $j=0,1,\cdots,k+1$.  
Elementary algebraic calculations then yield \eqref{coe}-\eqref{f-grow} for certain constants $L_1,L_2,L_3,L_4$ and $q=2k$, provided $a_{2k+1} + c_{k+1}^2 \sum_{m \in \nn_+} q_m \|g_m\|_{L_\xi^\infty}^2 < 0$.
In particular, this includes the double-well potential drift $f(\xi)=\epsilon^{-2}(\xi-\xi^3)$, $\xi \in \mathbb{R}$, $0<\epsilon < 1$, corresponding to the SACE \eqref{ac}, and a quadratic-growth diffusion function $g$ with a sufficiently small leading coefficient. 
\end{remark}

\begin{remark} \label{rk-ap-lin}
If $g$ is of linear growth (including the additive noise case), then the condition \eqref{coe} decomposes as
\begin{align} 
\xi f(\xi) \le L_1'-L_2' |\xi|^{q+2}, & \quad \xi \in \rr, \label{f-coe} \\  
 \|G(u)\|_{\LL_2^0}^2 \le L_5'+ L_6' \|u\|^2, & \quad u \in H, \label{s-grow}
\end{align}
for some positive constants $L_1'$, $L_2'$, $L_5', L_6'$.
\end{remark}

    \subsection{Tamed Schemes}

To construct nonlinearity-explicit schemes for Eq.~\eqref{see} when the coercive condition \eqref{coe} involves superlinear drift and diffusion coefficients, we employ a taming technique that simultaneously tames both the nonlinear drift and the diffusion. 

More precisely, we consider the following nonlinearity-explicit tamed EM (TEM) scheme applied to Eq.~\eqref{see}:
 \begin{align} \label{tame}
 	Z_n=Z_{n-1}+A Z_n \tau
 	+ F_\tau (Z_{n-1}) \tau 
 	+ G_\tau (Z_{n-1})\delta_{n-1} W; \quad Z_0=X_0,
 \end{align}
 where $\tau \in (0, 1]$ is the temporal step-size, $\delta_{n-1} W:=W(t_n) -W(t_{n-1})$, $n \in \nn_+$, and the tamed drift and diffusion functions are defined by
 \begin{align} \label{bs-tau}
 f_\tau (\xi):=\frac{f(\xi)}{(1+\tau|\xi|^{2q})^{1/2}}, \quad 
 g_\tau(\xi):=\frac{g(\xi)}{(1+\sqrt{\tau}|\xi|^q)^{1/2}}, \quad 
\xi \in \rr.
\end{align}   
We refer to \cite{LW26a} for a similar taming technique in the finite-dimensional setting. 

The TEM scheme \eqref{tame} is linearly implicit and therefore admits a unique pathwise solution. Moreover, $\{Z_n, \mathcal F_{t_n}\}_{n \in \mathbb{N}}$ is clearly a (time-homogeneous) Markov chain.

\section{Unique Ergodicity of Tamed Schemes}
 \label{sec3}

This section investigates the numerical ergodicity of the TEM scheme \eqref{tame} and its Galerkin-based full-discretizations of Eq.~\eqref{see}, exploring their Lyapunov structures and probabilistic regularity.
 
\subsection{Lyapunov Structure}

 The one-step Lyapunov estimate in this part indicates that the TEM scheme \eqref{tame} preserves a Lyapunov structure as the original one of Eq. \eqref{see}, even though it is not a Lyapunov functional of the classical EM scheme as indicated in \cite{BHJKLS19}. This shows that \eqref{tame} is unconditionally stable over an infinite time horizon, improving the finite-time unconditional stability result in \cite{HS23} for another type of tamed scheme in the 1D case. 

The main observation is the following coercive estimate about the tamed drift and diffusion coefficients given in \eqref{bs-tau}.
Hereafter, we omit the integration variable when an integration is present to simplify the notation.

\begin{lemma} 
Under conditions \eqref{coe}-\eqref{f-grow}, there exist positive constants $T_1$, $T_2$, and $\tau_{\max} < 1$ such that for any $\tau \in (0, \tau_{\max}]$, 
\begin{align} \label{coe-tau}
\<u, F_\tau(u) \> + \tau \|F_\tau(u)\|^2+\|G_\tau(u)\|_{\LL_2^0}^2 
& \le T_1 - T_2 \|u\|^2.
\end{align} 
\end{lemma}

\begin{proof}
From the representation \eqref{bs-tau} and conditions \eqref{coe}-\eqref{f-grow}, we have 
\begin{align*}
& \<u, F_\tau(u) \> + \tau \|F_\tau(u)\|^2+\|G_\tau(u)\|_{\LL_2^0}^2 \\
& =\int_\OOO \frac{u f(u)}{(1+\tau|u|^{2q})^{1/2}}
+  \frac{\tau |f(u)|^2}{1+\tau|u|^{2q}}
+ \frac{|g(u) |^2 \sum_{m \in \nn_+} q_m |g_m|^2}{(1+\sqrt \tau |u|^q)} {\rm d}\xi \\
& \le \int_\OOO \frac{2\tau (L_3^2+L_4^2 |u|^{2q+2})}{1+\tau|u|^{2q}}
+ \frac{u f(u)+|g(u)|^2  \sum_{m \in \nn_+} q_m |g_m|^2}{(1+\tau|u|^{2q})^{1/2}} {\rm d}\xi  \\ 
& \le \int_\OOO \frac{2\tau (L_3^2+L_4^2 |u|^{2q+2})}{1+\tau|u|^{2q}}
+ \frac{L_1 - L_2 |u|^{q+2}}{(1+\tau|u|^{2q})^{1/2}} {\rm d}\xi  \\ 
& \le (L_1+2 L_3^2 \tau) |\OOO| + \int_\OOO \frac{2 L_4^2 \tau |u|^{2q+2}}{1+\tau|u|^{2q}} - \frac{L_2 |u|^{q+2}}{(1+\tau|u|^{2q})^{1/2}} {\rm d}\xi. 
\end{align*} 
Using the following estimate derived in \cite[(2.18) in Lemma 2.1]{LW26a}:
 \begin{align*} 
\frac{2L_4^2 \tau \xi^{2 q+2}}{1+\tau \xi^{2q}}
- \frac{L_2 \xi^{q+2}}{(1+\tau \xi^{2q})^{1/2}}  
\le  -\frac{\beta \alpha^q}{\sqrt{1+\alpha^{2q}}} \xi^2
+ \frac{\beta \alpha^{q+2}}{\sqrt{1+\alpha^{2q}}},
\quad \xi \ge 0, 
 \end{align*} 
 for all $\beta \in (0, L_2)$, $\tau \in (0, (L_2-\beta)^2 /(8 L_4^4)]$, and $\alpha>0$, we conclude \eqref{coe-tau} with $\tau_{\max}=\min\{L_2^2 /(8 L_4^4), 1\}$ and certain positive constants $T_1, T_2$.
 \qed 
\end{proof}

We have the following Lyapunov estimate for the TEM scheme \eqref{tame}, where $\ee_n [\cdot]:=\ee[\cdot |\mathcal F_{t_n}]$ denotes the conditional expectation with respect to $\mathcal F_{t_n}$, $n \in \nn$; we always assume that $X_0$ is an $H$-valued, $\mathcal F_0$-measurable r.v.

 \begin{theorem} \label{tm-lya}
Assume that Assumption \ref{ap} holds and $\ee \|X_0\|^2 < \infty$.
 There exist positive constants $T_3$, $T_4$, and $\tau_{\max} < 1$ such that for any $\tau \in (0, \tau_{\max}]$, 
 	\begin{align} \label{lya}
\ee_{n-1} \|Z_n\|^2 + 2 \tau \ee_{n-1} \|\nabla Z_n\|^2 
 	\le  T_3 \tau + (1- T_4 \tau) \|Z_{n-1}\|^2, \quad n \in \nn_+.
 	\end{align}
 \end{theorem}

 \begin{proof}
 Let $n \in \nn_+$.
Testing \eqref{tame} with $Z_n$ and using the elementary equality 
\begin{align} \label{ab}
\<\alpha-\beta, \alpha\>
=(|\alpha|^2-|\beta|^2)/2+|\alpha-\beta|^2/2,
\quad \alpha, \beta\in \rr,
\end{align}  
and Cauchy--Schwarz inequality, we obtain  
\begin{align*}
& \|Z_n\|^2-\|Z_{n-1}\|^2 
+  \|Z_n-Z_{n-1}\|^2 +2 \tau \|\nabla Z_n\|^2 \nonumber  \\
&=2\tau \<Z_n-Z_{n-1}, F_\tau(Z_{n-1}) \> 
+ 2\tau \<Z_{n-1}, F_\tau(Z_{n-1})\> \nonumber  \\
& \quad + 2 \<Z_n-Z_{n-1}, G_\tau (Z_{n-1}) \delta_{n-1} W\>
+ 2 \<Z_{n-1}, G_\tau (Z_{n-1}) \delta_{n-1} W\> \nonumber \\
& \le 2 \tau^2 \|F_\tau(Z_{n-1})\|^2
+ 2\tau \<Z_{n-1}, F_\tau(Z_{n-1})\> + \|Z_n-Z_{n-1}\|^2 \nonumber  \\
& \quad + 2 \|G_\tau (Z_{n-1}) \delta_{n-1} W\|^2
+ 2 \<Z_{n-1}, G_\tau (Z_{n-1}) \delta_{n-1} W\>.
\end{align*}  
This shows 
\begin{align} \label{lya-err}
& \|Z_n\|^2+2 \tau \|\nabla Z_n\|^2   \nonumber \\
&\le \|Z_{n-1}\|^2 + 2\tau^2 \|F_\tau(Z_{n-1})\|^2
+ 2\tau \<Z_{n-1}, F_\tau(Z_{n-1}) \> \nonumber  \\
& \quad + 2 \|G_\tau (Z_{n-1}) \delta_{n-1} W\|^2
+ 2 \<Z_{n-1}, G_\tau (Z_{n-1}) \delta_{n-1} W\>.
\end{align}  
  
Note that $Z_{n-1} \sim \mathcal F_{t_{n-1}}$ and $\delta_{n-1} W$ is independent of $\mathcal F_{t_{n-1}}$. 
Then, taking the conditional expectation  $\ee_{n-1} [\cdot]$ on both sides of the above equation and using It\^o isometry and the fact 
$\ee_{n-1} \|G_\tau (Z_{n-1}) \delta_{n-1} W\|^2 
= \tau \|G_\tau (Z_{n-1})\|_{\LL_2^0}^2,$ 
we derive  
 \begin{align*} 
& \ee_{n-1} \|Z_n\|^2 + 2 \tau \ee_{n-1} \|\nabla Z_n\|^2  \\
&\le \|Z_{n-1}\|^2  
+ 2\tau [\<Z_{n-1}, F_\tau(Z_{n-1}) \> + \tau \|F_\tau(Z_{n-1})\|^2+\|G_\tau (Z_{n-1})\|_{\LL_2^0}^2].
\end{align*} 
We conclude \eqref{lya} by \eqref{coe-tau} with $T_3:=2 T_1$ and $T_4:=2 T_2$ for all $\tau \in (0, \tau_{\max}]$. 
\qed
\end{proof}

 \begin{remark}  
Recall that a Borel function $V: H \to [0, +\infty]$ is called a Lyapunov functional of a Markov chain $\{(Z_n, \mathcal F_{t_n})\}$, if  its sublevel sets
$K_R:= \{x \in H : V(x) \leq R\}$ are compact in $H$ for any sufficiently large $R>0$ and the drift condition 
$P V(x) := \mathbb E [V(Z_n) \mid Z_{n-1}=x] \le \beta V(x) + \kappa$ holds for some constants $\beta \in (0, 1)$ and $\kappa>0$ and for any $x \in H$, $n \in \mathbb N_+$. 
The drift condition yields that $a_n \le \beta a_{n-1} + \kappa$ with $a_n = \mathbb{E} V(Z_n)$, $n \in \mathbb N_+$.
  Recursively iterating the previous inequality, we obtain 
$a_n \le \beta^n a_0 + \kappa \sum_{k=0}^{n-1} \beta^k
        \le a_0 + \frac{\kappa}{1-\beta}$, $n \in \mathbb N_+$.
    Thus, taking the supremum over \(n\) shows that $\{Z_n^0\}$ (starting from $X_0=0$) is unconditionally stable on an infinite time horizon:
$\sup_{n \in \mathbb{N}_+} \mathbb{E} V(Z_n^0) \le V(0) + \frac{\kappa}{1-\beta}$. 
    In applications, this drift condition can be weakened as $\sup_{\mu \in \Lambda} \mu(V) < +\infty$, with $\mu(V):=\int_H V(x)\,\mu(dx)$, for a sequence $\Lambda$ of probability measures (usually taken as the empirical occupation measures $\mu_N(B):=\frac1N \sum_{n=1}^N \mathbb P(Z^0_n \in B)$, $N \in \mathbb N_+$, $B \in \mathcal B(H)$) in $H$, which yields the tightness of $\Lambda$ (see, e.g., \cite[Proposition 6.8]{Da06}).

The estimate \eqref{lya}, together with the compact embedding $\dot H^1 \subset H$, implies that $V(\cdot): = \|\cdot\|_1^2$ (or $V(\cdot) : = \|\cdot\|^2 + 2\tau \|\nabla \cdot\|^2$) is also a Lyapunov functional for the TEM scheme \eqref{tame}. 
Indeed, we take $\tau < \min\{\tau_{\max}, 1/T_4\}$ so that $\beta:=1-T_4 \tau \in (0,1)$ and $P V(x) \le \beta V(x) + \kappa$ with $V(\cdot) : = \|\cdot\|^2 + 2\tau \|\nabla \cdot\|^2$ and $\kappa=T_3 \tau >0$ by \eqref{lya}. 
    On the other hand, taking the expectation in \eqref{lya} and iterating shows that the scheme \eqref{tame} is unconditionally stable on an infinite time horizon:
\begin{align*}
E \|Z_n\|^2 & \le  T_3 \tau + (1- T_4 \tau) E \|Z_{n-1}\|^2
\le \cdots \\
& \le T_3 \tau [1+(1- T_4 \tau)+\cdots+(1- T_4 \tau)^{n-1}] 
+ (1- T_4 \tau)^n E \|Z_0\|^2 \\
& \le  T_3/T_4 + e^{-T_4 \tau n} E \|X_0\|^2, \quad \forall~ n \in \nn_+.
      \end{align*} 
Similarly, taking expectation in \eqref{lya} and summing over \(1 \le n \le N\) gives
\[
\sum_{n=1}^{N} \mathbb E\|Z_n\|^2 + 2\tau \sum_{n=1}^{N} \mathbb E\|\nabla Z_n\|^2
\le T_3\tau N + (1 - T_4\tau)\sum_{n=0}^{N-1} \mathbb E\|Z_n\|^2.
\] 
Rearranging and discarding certain terms, followed by division by $2N$, yields
\begin{align*}
\sup_{N \in \nn_+} \frac{1}{N}\sum_{n=1}^{N} \mathbb E\|\nabla Z_n\|^2 \tau 
& \le \frac{T_3+\mathbb E\|X_0\|^2}{2}.
\end{align*}  
Consequently, for $V(\cdot): = \|\cdot\|_1^2$, we obtain 
\begin{align*}
\sup_{N \in \nn_+} \mu_N(V)
& =\sup_{N \in \nn_+} \frac{1}{N}\sum_{n=1}^{N} E V(Z_n^0) \\
& \le  \sup_{N \in \nn_+} \mathbb E\|Z_N^0\|^2
+ \sup_{N \in \nn_+} \frac{1}{N}\sum_{n=1}^{N} \mathbb E\|\nabla Z_n^0\|^2  \\
& \le \frac{T_3}{T_4} + \frac{T_3}{2 \tau}<\infty,
\end{align*}  
thereby, in combination with the classical result \cite[Proposition 6.8]{Da06}, showing the tightness of $\{\mu_N\}_{N \in \mathbb{N}_+}$.
This is comparable to the Lyapunov estimate for Eq. \eqref{see} obtained, for instance, in \cite[Theorem 2.2]{HLL20}.  
\end{remark}

 \begin{remark}\label{rk-lya} 
 If the diffusion function $g$ is of linear growth, then we do not need to tame $g$ and consider the drift-TEM scheme
\begin{align} \label{tem-}
 	Y_n=Y_{n-1}+ A Y_n
	+ F_\tau (Y_{n-1}) \tau
 	+ G(Y_{n-1})\delta_{n-1} W; \quad Y_0=X_0,
 \end{align} 
 $n \in \nn_+$. 
Similar arguments used in the proof of Theorem \ref{tm-lya} yield  
 	\begin{align} \label{lya-}
\ee_{n-1} \|Y_n\|^2 + 2 \tau \ee_{n-1} \|\nabla Y_n\|^2
&\le T_3' \tau + (1- T_4' \tau) \|Y_{n-1}\|^2, \quad n \in \nn_+,
 	\end{align}   
for any $\tau \in (0, {\widetilde \tau}_{\max}]$ with certain ${\widetilde \tau}_{\max} \in (0, 1]$, and for some	positive constants $T_3'$ and $T_4'$, 
provided conditions \eqref{f-grow} and \eqref{f-coe}-\eqref{s-grow} hold. 
 \end{remark}

\subsection{Application to Galerkin-based tamed full-discretizations}

In this part, we extend the arguments and results of the previous part to construct tamed full discretizations using Galerkin-based spatial approximations, without imposing any additional conditions. Although our method can be generalized to spectral Galerkin approximations, we focus here on the Galerkin finite element method (FEM) applied to Eq.~\eqref{see}.

Let $h\in(0,1)$ and let $\mathcal{T}_h$ be a regular family of quasi-uniform partitions of $\mathcal{O}$ with maximum mesh size $h$. Denote by $V_h\subset V$ the space of continuous functions on $\overline{\mathcal{O}}$ that are piecewise linear with respect to $\mathcal{T}_h$ and vanish on $\partial\mathcal{O}$. Let $\mathcal{B}(V_h)$ be the Borel $\sigma$-algebra generated by $V_h$. Furthermore, let $A_h: V_h\to V_h$ and $\mathcal{P}_h: V^* \to V_h$ denote the discrete Dirichlet Laplacian and the generalized orthogonal projection operator, respectively, defined by 
\begin{align*}  
\<A_h u^h, v^h\> & =-\<\nabla u^h, \nabla v^h\>,
\quad u^h, v^h\in V_h,  \\
\<\PP_h u, v^h\> & =_1\<v^h, u\>_{-1},
\quad u\in V^*,\ v^h\in V_h. 
\end{align*}

We discretize the TEM scheme \eqref{tame} in space using FEM. That is, we seek a sequence $\{Z_n^h\}$ of $V_h$-valued Markov chains adapted to $\{\mathcal F_{t_n}\}$ such that 
 \begin{align} \label{tfem}
 	Z_n^h=Z_{n-1}^h+A_h Z_n^h \tau
 	+ \PP_h F_\tau (Z_{n-1}^h) \tau 
 	+ \PP_h G_\tau(Z_{n-1}^h)\delta_{n-1} W,   
 \end{align}
 $n \in \nn_+$, with initial datum $Z^h_0=\PP_h X_0$.
 We call this scheme the tamed-FEM. It can be equivalently written as   
\begin{align}\label{full}
Z_n^h=S_{h,\tau} Z_{n-1}^h+\tau S_{h,\tau} \PP_h F_\tau(Z_{n-1}^h)
+S_{h,\tau} \PP_h G_\tau(Z_{n-1}^h) \delta_{n-1} W, 
\end{align}
 $n \in \nn_+$, where $S_{h,\tau}:=({\rm Id}-\tau A_h)^{-1}$ is a space-time approximation of the continuous semigroup $S$ in one step.
Iterating \eqref{full} for $n-1$ times, we obtain  
\begin{align*}
Z^h_n 
=S_{h,\tau}^n  X^h_0+\tau \sum_{i=0}^{n-1}  S_{h,\tau}^{n-i} \PP_h F_\tau(Z^h_i)
+\sum_{i=0}^{n-1}  S_{h,\tau}^{n-i} \PP_h G_\tau (Z^h_i) \delta_i W,
\quad n \in \nn_+.
\end{align*} 
Similarly, in the linear growth diffusion case, we have the following equivalent drift-tamed-FEM schemes for $n \in \nn_+$, with $Y_0^h=\PP_h X_0$:
 \begin{align}  
Y_n^h & =Y_{n-1}^h + A_h Y_n^h \tau
+ \PP_h F_\tau (Y_{n-1}^h) \tau
 	+ \PP_h G(Y_{n-1}^h) \delta_{n-1} W,   \label{tfem-} \\
Y^h_n & =S_{h,\tau}^n  X^h_0+\tau \sum_{i=0}^{n-1}  S_{h,\tau}^{n-i} \PP_h F_\tau(Y^h_i)
+\sum_{i=0}^{n-1}  S_{h,\tau}^{n-i} \PP_h G(Y^h_i) \delta_i W. 	\label{full-sum}
 \end{align}

Similarly to the proof of Theorem \ref{tm-lya} and in light of Remark \ref{rk-lya}, we can establish the following Lyapunov estimates of the tamed-FEM \eqref{tfem} and the drift-tamed-FEM \eqref{tfem-}, respectively. The constants $T_3$, $T_4$, $\tau_{\max}$, and $T_3'$, $T_4'$, ${\widetilde \tau}_{\max}$ are those appearing in \eqref{lya} and \eqref{lya-}, respectively.

 \begin{theorem} \label{tm-lya-gtem}
Assume that Assumption \ref{ap} holds and $\ee \|X_0\|^2 < \infty$.
For any $\tau \in (0, \tau_{\max}]$ and $h \in (0, 1)$, the tamed-FEM \eqref{tfem} satisfies
 	\begin{align*}
\ee_{n-1} \|Z_n^h\|^2 + 2 \tau \ee_{n-1} \|\nabla Z_n^h\|^2 
 	\le T_3 \tau + (1- T_4 \tau) \|Z_{n-1}^h\|^2, \quad n \in \nn_+.  
 	\end{align*}   
 \end{theorem}

  \begin{theorem} \label{tm-lya-gtem-}
Assume that conditions \eqref{f-grow} and \eqref{f-coe}-\eqref{s-grow} hold and $\ee \|X_0\|^2 < \infty$.
For any $\tau \in (0, {\widetilde \tau}_{\max}]$ and $h \in (0, 1)$, the drift-tamed-FEM \eqref{tfem-} satisfies
 	\begin{align*} 
	\ee_{n-1} \|Y_n^h\|^2 + 2 \tau \ee_{n-1} \|\nabla Y_n^h\|^2
	 	\le T_3' \tau + (1- T_4' \tau) \|Y_{n-1}^h\|^2, \quad n \in \nn_+. 
 	\end{align*} 
 \end{theorem}

\subsection{Unique Ergodicity of Tamed-FEM and Drift-Tamed-FEM}

In this part, we focus on the probabilistic regularity and unique ergodicity of the tamed-FEM \eqref{tfem} and the drift-tamed-FEM \eqref{tfem-} under the following additional nondegenerate condition.
Consequently, our numerical ergodicity analysis does not apply to the standard linear multiplicative noise case, which is degenerate.

\begin{assumption} \label{ap-non} 
For any $u \in H$, $G(u) G(u)^\tr$ is a positive-definite linear operator in $H$, i.e.,
$\< G(u) G(u)^\tr v, v\> >0$, $v \in H$.
\end{assumption}
Recall that at the end of Section~\ref{sec2.1} we assumed that $Q$ is self-adjoint and positive-definite. It is well known that $Q^{1/2}$ is then also self-adjoint and positive-definite. One can easily show the equivalence between the positive-definiteness of $G(u)G(u)^\top$ and that of $G(u)Q^{1/2}[G(u)Q^{1/2}]^\top$; see \cite[Remark 4]{LL24}.

As noted earlier in this section, $\{Z_n^h, \mathcal F_{t_n}\}_{n\in\mathbb{N}}$ and $\{Y_n^h, \mathcal F_{t_n}\}_{n\in\mathbb{N}}$ generated by the tamed-FEM \eqref{tfem} and the drift-tamed-FEM \eqref{tfem-}, respectively, are Markov chains.  
Let $P: V_h \times \mathcal{B}(V_h) \to [0,1]$ and $\widetilde{P}: V_h \times \mathcal{B}(V_h) \to [0,1]$ denote the transition kernels of $\{Z_n^h, \mathcal F_{t_n}\}_{n\in\mathbb{N}}$ and $\{Y_n^h, \mathcal F_{t_n}\}_{n\in\mathbb{N}}$, respectively. Then for any $n \in \nn$, $x \in V_h$, and $A \in \BB(V_h)$,
 \begin{align} \label{kernel} 
 P(x,A)=\mathbb{P}(Z_{n+1}^h \in A~|~Z_n^h=x), ~
{\widetilde P}(x,A)=\mathbb{P}(Y_{n+1}^h \in A~|~Y_n^h=x). 
 \end{align}
 We refer to \cite{LL24} for the definitions and properties of invariant measure and ergodicity for Markov chains.

We have the following probabilistic regularity for the transition kernels $P$ and ${\widetilde P}$ defined in \eqref{kernel} associated with \eqref{tfem} and \eqref{tfem-}, respectively.
 
 \begin{proposition} \label{prop-reg}
 Let Assumption \ref{ap-non} hold. For any $\tau \in (0,1)$, the transition kernels $P$ and ${\widetilde P}$ are regular on $\BB(V_h)$, i.e., for every $x \in V_h$, the measures $P(x,\cdot)$ and $\widetilde{P}(x,\cdot)$ are each equivalent to a common reference measure. Consequently, each of the Markov chains --- the tamed-FEM \eqref{tfem} and the drift-tamed-FEM \eqref{tfem-} --- admits at most one invariant measure on $\mathcal{B}(V_h)$.
 \end{proposition}
 
 \begin{proof}
For $ x \in V_h$ and $A \in \BB(V_h)$, from the tamed-FEM \eqref{tfem} we have \begin{align*}
P(x,A)&=\mathbb{P} (Z_{n+1}^h \in A |Z_n^h=x)\\
&=\mathbb{P} (x+ \PP_h F_\tau(x)
+\PP_h G_\tau(x) \delta_n W) \in ({\rm Id}- \tau A_h) (A) )\\
&=\mu_{x+\PP_h F_\tau(x),  [\PP_h G_\tau(x)] [\PP_h G_\tau(x)]^* \tau}(({\rm Id}- \tau A_h) (A)),
 \end{align*}
where the last equality uses that $x+\PP_h F_\tau(x)+ \PP_h G_\tau(x) \delta_n W$ is normally distributed with mean $x+\PP_h F_\tau(x)$ and covariance operator $[\PP_h G_\tau(x)] [\PP_h G_\tau(x)]^* \tau$, and $\mu_{z, C}$ denotes the Gaussian measure on $V_h$ with mean $z$ and covariance operator $C$. 

Under Assumption~\ref{ap-non}, the operator $[\mathcal{P}_h G_\tau(x)][\mathcal{P}_h G_\tau(x)]^* \tau$ is nondegenerate in $\mathcal{L}(V_h)$ (see \cite[Remark 4]{LL24}). Consequently, by Feldman--H\'ajek theorem on the finite-dimensional space $V_h$, all Gaussian measures of the form $\{\mu_{x, [\PP_h G(x) {Q^{1/2}}] [\PP_h G(x) {Q^{1/2}}]^* \tau}: x\in V_h\}$ are equivalent. This shows that $P$ is regular. The regularity of ${\widetilde P}$ follows analogously, which completes the proof.
\qed
\end{proof}

Combining the Lyapunov structure from Theorem~\ref{tm-lya} and Remark~\ref{rk-lya} with the probabilistic regularity established in Proposition~\ref{prop-reg}, we obtain the following ergodicity results for the tamed-FEM \eqref{tfem} and the drift-tamed-FEM \eqref{tfem-}, respectively.
We also note that, following the arguments in \cite[Theorem 2]{LL24} and \cite[Theorem 2.5]{MSH02}, one could prove geometric ergodicity under additional differentiability assumptions on \(f\) and \(g\) by verifying a minorization condition.

\begin{theorem}  
Let Assumptions \ref{ap} and \ref{ap-non} hold.  
The tamed-FEM \eqref{tfem} is uniquely ergodic in $\BB(V_h)$ for any $\tau \in (0, \tau_{\max}]$.   

\end{theorem}

\begin{theorem}  
Let conditions \eqref{f-grow} and \eqref{f-coe}-\eqref{s-grow} and Assumption \ref{ap-non} hold.   
The drift-tamed-FEM \eqref{tfem-} is uniquely ergodic in $\BB(V_h)$ for any $\tau \in (0, {\widetilde \tau}_{\max}]$.  

\end{theorem}

 \section{Strong Error Estimates in Multiplicative Noise Case}
 \label{sec4}

In this section, we carry out a quantitative error analysis between the drift-tamed-FEM \eqref{tfem-} and Eq.~\eqref{see} under the $L^{2p}(\Omega)$-norm for $p\ge 2$, to demonstrate the effect of these tamed methods on the dynamics of \eqref{see-fg}.

To establish the desired strong error estimates, we impose the following conditions on the coefficients of Eq.~\eqref{see}, which are consistent with those in Assumptions~\ref{ap} and~\ref{ap-non}. We refer to \cite[Remarks 2.2-2.3 and Example 5.1]{LQ21} for examples under which all four assumptions are satisfied.

 \begin{assumption} \label{ap-f}
 $f \in \CC^1(\rr; \rr)$ and there exist positive constants $K_1, K_2$ such that for all $\xi \in \rr$ and $\tau \in (0, 1)$,  
\begin{align} 
|f'(\xi)| & \le K_1 (1+|\xi|^q), \label{f'} \\ 
[1+\tau |\xi|^{2q}]f'(\xi) - q \tau |\xi|^{2(q-1)} \xi f(\xi) 
& \le K_2 (1+\tau|\xi|^{2q})^{3/2}. \label{con-f'up} 
\end{align} 
\end{assumption} 

The duality method and the condition \eqref{f'} imply that (see \cite[(4.16)]{LQ21}) 
 \begin{align}  \label{F-}
\|F(u)-F(v)\|_{-1} & \le C (1+\|u\|^q_1+\|v\|^q_1 ) \|u-v\|, 
\quad u, v \in \dot H^1.   
\end{align}

\begin{remark} \label{rk-upper}
The condition \eqref{con-f'up} is equivalent to the uniform bound $f_\tau' \le K_2$, which is crucial for deriving optimal strong error estimates in this and the next sections.
Under this condition, letting $\tau \to 0$ also yields $f' \le K_2$, so that there exists a constant $\widetilde{K}_2$ such that 
\begin{align} \label{f-mon} 
(f(\xi)-f(\eta)) (\xi-\eta) & \le K_2 |\xi-\eta|^2, \quad 
\xi f(\xi) \le \widetilde{K}_2 (1+ \xi^2), \quad \xi, \eta \in \rr.  
\end{align} 
Moreover, \eqref{con-f'up} includes all odd polynomials with negative leading coefficients, as described in Remark~\ref{rk-ap-pol} (i.e., \eqref{bs-ex} with $a_{2k+1}<0$, $k \in \mathbb{N}_+$). In this case, direct calculations yield
$$D(\xi):=[1+\tau |\xi|^{2q}]f'(\xi)-q \tau |\xi|^{2(q-1)} \xi f(\xi)=f'(\xi)+\tau H(\xi), \quad \xi \in \rr,$$
where $f'$ and 
$H(\xi):=|\xi|^{2q} f'(\xi)-q |\xi|^{2(q-1)} \xi f(\xi)$, $\xi \in \rr$,
are polynomials of degrees $2k$ and $3q = 6k$, respectively, with leading coefficients $(2k+1)a_{2k+1}<0$ and $a_{2k+1}<0$, respectively. Consequently, both $f'$ and $H$ are bounded above, and hence $D$ admits a uniform upper bound, so that \eqref{con-f'up} holds.
\end{remark}

\begin{assumption} \label{ap-g}
\begin{enumerate}
\item[(1)]
$G: H\rightarrow \LL_2^0$ is Lipschitz continuous and $G: \dot H^1 \rightarrow \LL_2^1$ grows linearly, i.e., there exist positive constants $K_3, K_4$ such that 
\begin{align}   
\|G(u)-G(v)\|_{\LL_2^0} \le K_3 \|u-v\|, & \quad u,v \in H, \label{g-lip} \\ 
\|G(z) \|_{\LL_2^1} \le K_4 (1 + \|z\|_1), & \quad z \in \dot H^1. \label{g1} 
\end{align} 

\item[(2)] 
There exist positive constants $\theta$ and $K_5$ such that $G(\dot H^{1+\theta}_x)\subset \LL_2^{1+\theta}$ and 
\begin{align} \label{g2}
\|G(z) \|_{\LL_2^{1+\theta}}
\le K_5 (1+ \|z\|_{1+\theta}),
\quad z \in \dot H^{1+\theta}.
\end{align}
\end{enumerate}
\end{assumption}

It is known that under Assumption \ref{ap-g}(1),  Eq. \eqref{see} with $X_0 \in \dot H^1$ admits a unique mild solution in $L^p(\Omega; \CC([0, T]; \dot H^1)$ for arbitrary $p \ge 1$ that satisfy (see \cite{LQ21})
\begin{align} \label{mild}
X(t) &=S(t) X_0+\int_0^t S(t-r) F(u) {\rm d}r
+\int_0^t S(t-r) G(u) {\rm d}W(r),
\end{align}
$t \ge 0$. 
Moreover, we note that under Assumptions \ref{ap}, \ref{ap-non}, \ref{ap-f}, \ref{ap-g} and certain boundedness conditions on the pseudo-inverse of $G$, Eq. \eqref{see} is uniquely ergodic in $H$; see, e.g., \cite[Theorem 2.2]{HLL20}.

We will need the following ultracontractive and smoothing properties of the analytic $\CC_0$-semigroup $S$ (see, e.g., \cite[Theorem 6.13 in Chapter 2]{Paz83}):
\begin{align}
\|S(t) u\|_\mu \le C t^{-\frac{\mu-\nu}2} \|u\|_\nu,   
& \quad \forall~ t > 0, ~ 0 \le \nu \le \mu \le 2,~ u \in \dot H^\nu, \label{ana} \\
\|(S(t)-{\rm Id}_H) u\| \le C t^{\frac\rho2} \|u\|_\rho, \label{ana-hol} 
& \quad \forall~ t > 0, ~ 0 \le \rho \le 2, ~ u \in \dot H^\rho,  \\
\|S(t) u\|_{L_\xi^p} \le C \|u\|_{L_\xi^p},   
& \quad \forall~ t > 0, ~ 1 \le p \le \infty,~ u \in L_\xi^p. \label{ana-lp}
\end{align}   
Hereafter, we denote by $C$ a generic constant independent of the discretization parameters, which may take different values from one occurrence to another.
Similarly to \eqref{ana} (with $\mu=\nu=0$), there holds that  
\begin{align} \label{sht1}
\|S_{h,\tau}^k \PP_h u\|
\le \|u\|,
& \quad \forall~ k \in \nn_+, ~ u \in H. 
\end{align}

\subsection{A Priori Estimates}

In the following estimates, we will frequently use the properties of the tamed function $f_\tau$ defined in \eqref{bs-tau} and its associated Nemytskii operator $F_\tau$.

\begin{lemma} \label{lm-ftau} 
Under Assumption \ref{ap-f}, there exists a constant $C$ such that for any $\xi \in \rr$, $u, v \in \dot H^1$, and $z \in L_\xi^{2(3q+1)}$,
\begin{align} 
\<F_\tau(u), u\>_1 & \le C (1+\|u\|_1^2),  \label{Ftau-coe1}  \\ 
_{-1}\<F_\tau(u)-F_\tau(v), u-v\>_1 
& \le C \|u-v\|^2, \label{Ftau-mon} \\
|f_\tau'(\xi)| & \le  C \tau^{-1/2}, \quad \xi \in \rr, \label{ftau'-grow}  \\  
\tau \|F_\tau(u)\|_1^2
& \le C(1+\|u\|_1^2), \quad u \in \dot H^1, \label{tauFtau1}  \\
\|F_\tau(u)-F_\tau(v)\|_{-1} 
& \le C (1+\|u\|^q_1+\|v\|^q_1 ) \|u-v\|, 
\quad u, v \in \dot H^1, \label{Ftau-} \\ 
\|F(u)-F_\tau(u)\|_{-1} 
& \le C \tau^{1/2} (1+\|u\|_1^{2q+1}), \quad u \in \dot H^1, \label{F-Ftau} \\
\|F(z)-F_\tau(z)\|
& \le C \tau (1+\|z\|_{L_\xi^{2(3q+1)}}^{3q+1}), \quad z \in L_\xi^{2(3q+1)}. \label{F-Ftau+} 
\end{align}   
\end{lemma}

\begin{proof}   
 Let $\xi \in \rr$, $u, v \in \dot H^1$, and $z \in L_\xi^{2(3q+1)}$. 
When \eqref{con-f'up} holds, by Remark \ref{rk-upper}, $f_\tau' \le K_2$.
Consequently, we derive \eqref{Ftau-coe1} and \eqref{Ftau-mon}, respectively:
\begin{align*} 
\<F_\tau(u), u\>_1
& =\<F_\tau(u), u\>+\<F_\tau'(u) \nabla u, \nabla u\> \\
& =\int_D \Big(\frac{u f(u)}{(1+\tau|u|^{2q})^{1/2}}+f_\tau'(u) |\nabla u|^2 \Big){\rm d}\xi   \\
 & \le \widetilde{K}_2 (1+\|u\|^2) + K_2 \|\nabla u\|^2 \\
& \le C (1+\|u\|_1^2), \\
_{-1}\<F_\tau(u)-F_\tau(v), u-v\>_1 
& =\int_D f_\tau'(u(\xi)+\lambda (v(\xi)-u(\xi)))) |u(\xi)-v(\xi)|^2 {\rm d}\xi \\
& \le K_2 \|u-v\|^2, \quad \text{with some~} \lambda  \in (0, 1).
\end{align*} 
 
Direct calculations yield 
\begin{align*}
f_\tau'(\xi)=\frac{[1+\tau |\xi|^{2q}]f'(\xi)-q \tau |\xi|^{2(q-1)} \xi f(\xi)} 
{(1+\tau|\xi|^{2q})^{3/2}}.
\end{align*}
Then the condition \eqref{f'} and Young inequality imply \eqref{ftau'-grow}: 
\begin{align*}
|f_\tau'(\xi)|
& \le \frac{K_1 (1+\tau |\xi|^{2q}) (1+|\xi|^q)} 
{(1+\tau|\xi|^{2q})^{3/2}}
+ \frac{C \tau (1+|\xi|^{3q})} 
{(1+\tau|\xi|^{2q})^{3/2}}
\le C \tau^{-1/2}.
\end{align*}
This bound, in combination with the fact $\tau |f_\tau(\xi)|^2 \le C(1+|\xi|^2)$, yields \eqref{tauFtau1}:
\begin{align*} 
\tau \|F_\tau(u)\|_1^2
&=\tau \|F_\tau(u)\|^2+\tau \|F_\tau'(u) \nabla u\|^2 
\le C(1+\|u\|_1^2). 
\end{align*}

To show \eqref{Ftau-}-\eqref{F-Ftau}, we use the representation \eqref{bs-tau} to get 
\begin{align*}
f_\tau(\xi)-f_\tau(\eta)  
&= \frac{f(\xi) (1+\tau|\eta|^{2q})^{1/2}-f(\eta)(1+\tau|\xi|^{2q})^{1/2}}{(1+\tau|\xi|^{2q})^{1/2} (1+\tau|\eta|^{2q})^{1/2}} \nonumber  \\
&= \frac{f(\xi)-f(\eta)}{(1+\tau|\xi|^{2q})^{1/2}} 
+ \frac{f(\eta) [(1+\tau|\eta|^{2q})^{1/2}-(1+\tau|\xi|^{2q})^{1/2}]}{(1+\tau|\xi|^{2q})^{1/2} (1+\tau|\eta|^{2q})^{1/2}}  \nonumber   \\
&= \frac{f(\xi)-f(\eta)}{(1+\tau|\xi|^{2q})^{1/2}} 
+\tau f(\xi, \eta),
\end{align*}
where  
\begin{align*}
f(\xi, \eta):=\frac{f(\eta) (|\eta|^{2q}-|\xi|^{2q})} 
{(1+\tau|\xi|^{2q})^{1/2} (1+\tau|\eta|^{2q})^{1/2}[(1+\tau|\eta|^{2q})^{1/2}+(1+\tau|\xi|^{2q})^{1/2}]},
\end{align*}
$\eta \in \rr$.
As $f$ satisfies \eqref{f'}, we have 
\begin{align*}
|f_\tau(\xi)-f_\tau(\eta)|  
& \le \frac{|f(\xi)-f(\eta)|}{(1+\tau|\xi|^{2q})^{1/2}} + \frac{\tau |f(\eta)| \cdot ||\eta|^{2q}-|\xi|^{2q}|} 
{1+\tau|\xi|^{2q}+\tau|\eta|^{2q}}  \\  
& \le |f(\xi)-f(\eta)|
+ \frac{C \tau (1+|\eta|^{q+1}) \cdot |\xi-\eta| (1+|\eta|^{2q-1}+|\xi|^{2q-1}|} 
{1+\tau|\xi|^{2q}+\tau|\eta|^{2q}}  \\  
& \le |f(\xi)-f(\eta)|
+ \frac{C \tau (1+ |\xi|^{3q}+ |\eta|^{3q})}{1+\tau|\xi|^{2q}+ \tau |\eta|^{2q}} \cdot |\xi-\eta| \\ 
&\le C |\xi-\eta| (1+ |\eta|^q+ |\xi|^q),
\end{align*}
where we used the elementary inequality
$||\eta|^{2q}-|\xi|^{2q}| \le C |\xi-\eta| (1+|\eta|^{2q-1}+|\xi|^{2q-1})$ for some $C=C(q)>0$ and any $\xi, \eta \in \mathbb R$.
By the embedding \eqref{emb}, we get \eqref{Ftau-}:
\begin{align*} 
\|F_\tau(u)-F_\tau(v)\|_{-1}
& \le C \|F_\tau(u)-F_\tau(v)\|_{L_\xi^{[2(q+1)]'}}
\le C \|u-v\| (1+\|u\|_1^q+\|u\|_1^q).
\end{align*}  

Finally, it was shown in \cite[Remark 5]{LW26a} that 
\begin{align} \label{f-ftau}
|f(\xi)-f_\tau(\xi)|
& =\frac{\tau |\xi|^{2q} |f(\xi)|}{(1+\tau |\xi|^{2q})^{1/2} [(1+\tau |\xi|^{2q})^{1/2}+1]} 
 \le \frac12 \tau |\xi|^{2q} |f(\xi)|,
\end{align}
so that \eqref{F-Ftau} holds:
\begin{align*} 
\|F(z)-F_\tau(z)\|
 & \le C \tau \|1+|z|^{3q+1}\|_{L_\xi^2}
\le C \tau (1+\|z\|_{L_\xi^{2(3q+1)}}^{3q+1}).
\end{align*}  
Moreover, the equality \eqref{f-ftau} also yields that 
$|f(\xi)-f_\tau(\xi)| 
\le \frac{\tau |\xi|^{2q} |f(\xi)|}{2 \tau^{1/2} |\xi|^q} 
\le \frac12 \tau^{1/2}|\xi|^q|f(\xi)|$ for $\xi \neq 0$.
It is clear that $f(\xi)=f_\tau(\xi)|$ when $\xi=0$, and thus  
\begin{align*}
|f(\xi)-f_\tau(\xi)| 
\le \frac12 \tau^{1/2}|\xi|^q|f(\xi)|.
\end{align*}  
Then, from the embedding \eqref{emb} and the above inequality, we have 
\begin{align*} 
\|F(z)-F_\tau(z)\|_{-1}  
\le C \tau^{1/2} \|1+|z|^{2q+1}\|_{L_\xi^{[2(q+1)]'}}
\le C \tau^{1/2} (1+\|z\|_1^{2q+1}).
\end{align*} 
This shows \eqref{F-Ftau}, thereby completing the proof.  
\qed
\end{proof}

We now recall the following regularity properties of the exact solution $X$ to Eq.~\eqref{see}. In addition, an $L_\xi^{2(3q+1)}$-estimate of $X$ is required to estimate the term $\|F(X)-F_\tau(X)\|_{L_\omega^p L_\xi^2}$ in \eqref{j23-}, using the bound \eqref{F-Ftau+}, for temporal higher order of strong convergence analysis in the additive noise case.

Hereafter, we denote by $\mathcal{C}_t^\alpha L_\omega^p \dot H^\beta$ the space of adapted stochastic processes $Z$ for which the seminorm $\|\cdot\|_{\mathcal{C}_t^\alpha L_\omega^p \dot H^\beta}$ is finite: 
\begin{align*}   
\|Z\|_{\CC_t^\alpha L_\omega^p \dot H^\beta}
:= \sup_{0\le s < t \le T} \frac{\|Z(t)-Z(s)\|_{L_\omega^p \dot H^\beta}} 
{|t-s|^\alpha} 
& <\infty.
\end{align*}

\begin{lemma} \label{lm-u}
\begin{enumerate}
\item[(1)]
Let $X_0 \in \dot H^{1+\gamma}$ with $\gamma\in [0,1]$ and Assumptions \ref{ap-f}-\ref{ap-g}(1) (and (2) if $\gamma=1$) hold.  
For any $p \ge 2$ and $\beta\in [0,1]$, there exists a constant $C$ such that 
\begin{align}  
\|X\|_{L_\omega^p L_t^\infty \dot H^{1+\gamma}}
& \le 
\begin{cases}
C e^{CT} (1+\|X_0\|_1), \quad & \gamma=0; \\
C e^{CT} (1+\|X_0\|^{q+1}_{1+\gamma}), \quad & \gamma\in [0,1); \\
C e^{CT} (1+\|X_0\|^{(q+1)^2}_2), \quad & \gamma=1,
\end{cases}  \label{reg-u} \\
\|X\|_{\CC_t^{\frac{1+\gamma-\beta}2\wedge \frac12} L_\omega^p \dot H^\beta}
& \le 
\begin{cases}
C e^{CT} (1+\|X_0\|^{q+1}_{1+\gamma}), \quad & \gamma\in [0,1); \\
C e^{CT} (1+\|X_0\|^{(q+1) ^2}_2), \quad & \gamma=1.
\end{cases}  \label{hol-u}
\end{align}  

\item[(2)]
Let $X_0 \in \dot H^1 \cap L_\xi^{2(3q+1)}$ and Assumptions \ref{ap-f}-\ref{ap-g}(1) hold. 
For any $p \ge 2$, there exists a constant $C$ such that  
\begin{align}   \label{reg-u+}
\|X\|_{L_\omega^p L_t^\infty L_\xi^{2(3q+1)}}
& \le  C e^{CT} (1+\|X_0\|_1+\|X_0\|_{L_\xi^{2(3q+1)}}^{q+1}).
\end{align}
\end{enumerate}  
\end{lemma}

\begin{proof}
(1) The moment estimate \eqref{reg-u} and the H\"older estimate \eqref{hol-u} were shown in Theorem 3.1+Proposition 3.1+Theorem 3.3 and Corollary 3.2 of \cite{LQ21}, respectively.

(2) It was shown in \cite[Proposition 3.1]{LQ21} that 
\begin{align*}
& \Big\|\int_0^t S(t-r)  F(X(r)) {\rm d}r
+\int_0^t S(t-r)  G(X(r)) {\rm d}W(r) \Big\|_{L_\omega^p L_t^\infty \dot H^{1+\beta}} \\
& \le C (\|F(X)\|_{L_\omega^p L_t^\infty L_\xi^2}
+ \|G(X)\|_{L_\omega^p L_t^\infty \LL_2^1}) \\
& \le C e^{CT} (1+\|X_0\|^{q+1}_1), \quad \forall~ \beta \in [0,1 ).
\end{align*}
This inequality, in combination with the embedding $\dot H^{1+\beta} \subset L^\infty$ for sufficiently large $\beta \in (1/2, 1)$ and the estimate on $J^0$ that  
\begin{align*} 
\|S(t)  X_0\|_{L_\omega^p L_t^\infty L_\xi^{2(3q+1)}}
\le C \|X_0\|_{L_\xi^{2(3q+1)}}, 
\end{align*} 
followed by \eqref{ana-lp} with $p=2(3q+1)$, completes the proof of \eqref{reg-u+}.  
\qed
\end{proof}

\begin{remark}
One might expect the estimates in Lemma~\ref{lm-u} to be uniform with respect to $T$ when the dissipative condition in Assumption~\ref{ap} is imposed. However, under the standard monotone/Lipschitz conditions \eqref{f-mon}/\eqref{g-lip}, the $T$-dependence of the strong convergence rates seems unavoidable in the general setting; we refer to \cite{Liu26} for uniform bounds and strong error estimates under a strong dissipativity condition. 
\end{remark}

\subsection{Auxiliary Process and Moments' Estimates}

As in \cite{LQ21}, we define the auxiliary process $\{{\widehat Y}_n^h\}$ by ${\widehat Y}_0^h=\PP_h X_0$ and  
\begin{align}\label{aux}
{\widehat Y}_n^h 
=S_{h,\tau}^n \PP_h X_0
+\tau \sum_{i=0}^{n-1} S_{h,\tau}^{n-i} \PP_h F_\tau(X(t_i))
+\sum_{i=0}^{n-1} S_{h,\tau}^{n-i} \PP_h G(X(t_i)) \delta_i W,
\end{align}
for $n \in \nn_+$, where the terms $Y_i^h$ in the discrete deterministic and stochastic convolutions of \eqref{full-sum} are both replaced by $X(t_i)$.
It is clear that for $n \in \nn_+$,
\begin{align} \label{aux+} 
{\widehat Y}_n^h
&={\widehat Y}_{n-1}^h+\tau A_h {\widehat Y}_n^h 
+\tau \PP_h F_\tau(X(t_{n-1}))
+\PP_h G(X(t_{n-1})) \delta_{n-1} W.
\end{align}

By further developing the idea from \cite[Proposition 3.1]{Liu22}, we obtain the following general moment estimate for the auxiliary process \eqref{aux}. This moment estimate in the $\dot H^1$-norm is crucial for deriving the higher-order moment estimates required for the error analysis in the next two sections.

\begin{lemma} \label{lm-reg-aux} 
Let $X_0 \in \dot H^1$ and Assumptions \ref{ap-f}+\ref{ap-g}(1) hold.  
For all $p \ge 2$, there exists a constant $C$ such that 
\begin{align}  \label{reg-aux} 
 \sup_{1 \le n \le N} \ee \|{\widehat Y}_n^h\|_1^{2p}  
& \le C e^{CT} (1+\|X_0\|_1^{2p(q+1)}).
\end{align}   
\end{lemma}

\begin{proof} 
 Let $p \ge 2$ and $n \in \nn_+$.
 Testing \eqref{aux+} with ${\widehat Y}_n^h$ under the $\dot H^1$-norm and using the elementary equality \eqref{ab}, Cauchy--Schwarz inequality, integration by parts formula, and the fact that $P_h$ is bounded in $H$ and $\dot H^1$, we obtain  
\begin{align*}
& \|{\widehat Y}_n^h\|_1^2-\|{\widehat Y}_{n-1}^h\|_1^2 
+  \|{\widehat Y}_n^h-{\widehat Y}_{n-1}^h\|_1^2 
+2 \tau \|\nabla {\widehat Y}_n^h\|_1^2 \nonumber  \\
&=2\tau \<{\widehat Y}_n^h, \PP_h F_\tau(X(t_{n-1}))\>_1
+ 2 \<{\widehat Y}_n^h-{\widehat Y}_{n-1}^h, \PP_h G(X(t_{n-1})) \delta_{n-1} W\>_1
+ M_{n-1}^h \\
& \le \tau \|\nabla {\widehat Y}_n^h\|_1^2 + C \tau \|F_\tau(X(t_{n-1}))\|^2 \\
& \quad + \|{\widehat Y}_n^h-{\widehat Y}_{n-1}^h\|_1^2
+ \|G (X(t_{n-1})) \delta_{n-1} W\|_1^2
+ M_{n-1}^h,
\end{align*}   
where $M_{n-1}^h:=2 \<{\widehat Y}_{n-1}^h, \PP_h G(X(t_{n-1})) \delta_{n-1} W\>_1$. 
It follows that 
\begin{align*}  
& \|{\widehat Y}_n^h\|_1^2-\|{\widehat Y}_{n-1}^h\|_1^2 + \tau \|\nabla {\widehat Y}_n^h\|_1^2 \\
& \le C \tau \|F_\tau(X(t_{n-1}))\|^2 
+ \|G (X(t_{n-1})) \delta_{n-1} W\|_1^2
+ M_{n-1}^h.
\end{align*}
Summing up the above inequality yields that  
\begin{align}  \label{lp-}
& \|{\widehat Y}_n^h\|_1^2 + \tau \sum_{k=1}^n \|\nabla {\widehat Y}_k^h \|_1^2 - \|X_0\|_1^2 \nonumber \\
& \le C \tau \sum_{k=1}^n \|F_\tau(X(t_{k-1}))\|^2
+ \sum_{k=1}^n \|G(X(t_{k-1})) \delta_{k-1} W\|_1^2
+ \sum_{k=1}^n M_{k-1}^h.
\end{align}
Then, we use the elementary inequality 
\begin{align} \label{in-sum}
\Big(\sum_{i=0}^{N-1} a_i \Big)^\gamma 
\le N^{\gamma-1} \sum_{i=0}^{N-1} a_i^\gamma,
\quad \forall~\gamma \ge 1, ~ a_i \ge 0,
\end{align} 
and the growth condition \eqref{f-grow} to derive 
\begin{align*}  
\ee\Big[C \tau \sum_{k=1}^N \|F_\tau(X(t_{k-1}))\|^2 \Big]^p  
& \le C \sum_{k=1}^N \ee \|F_\tau(X(t_{k-1}))\|^{2p} \tau \\
& \le C \sum_{k=1}^N (1+\ee \|X(t_{k-1})\|_1^{2p(q+1)}) \tau.
\end{align*}
Similarly, using Burkholder inequality in \cite[Theorem 6.3.6]{Str10}: for any $p >1$, there exist constants $C_1, C_2$ such that for any real-valued martingale $\{B_n, \mathcal{F}_{t_n}\}$,
\begin{align*}
C_1 \sup_{1 \le n \le N} \|B_n - B_0\|^p_{L^p_\omega}
\leq \ee \Big\| \sum_{n=1}^N |B_n - B_{n-1}|^2 \Big\|^{p/2}
\leq C_2 \sup_{1 \le n \le N} \|B_n - B_0\|^p_{L^p_\omega},
\end{align*}
with $B_n := \sum_{1 \le k \le n} M_{n-1}^h$, $n \in \mathbb{N}_+$, $B_0 := 0$, 
Minkowski, Cauchy--Schwarz, and Young inequalities, and the condition \eqref{g1}, we obtain  
\begin{align}  \label{est-sto}
& \ee\Big[\sum_{k=1}^N \|G(X(t_{k-1})) \delta_{k-1} W\|_1^2 \Big]^p  
+ \ee \sup_{1 \le n \le N} \Big|\sum_{k=1}^n M_{k-1}^h\Big|^p \nonumber \\
& \le C N^{p-1} \sum_{k=1}^N \ee \|G(X(t_{k-1})) \delta_{k-1} W\|_1^{2 p}  
+ C N^{\frac p2-1} \sum_{k=1}^N \|M_{k-1}^h\|_{L_\omega^p}^p \nonumber \\
& \le C \sum_{k=1}^N (1+\ee \|{\widehat Y}_{k-1}^h\|_1^{2p}+\ee \|X(t_{k-1})\|_1^{2p}) \tau.
\end{align} 

Taking $L^p_\omega$-norm on both sides of \eqref{lp-} and using the above two estimates, we derive
\begin{align*}  
& \ee \sup_{1 \le n \le N} \|{\widehat Y}_n^h\|_1^{2p}
+ \ee \Big(\sum_{k=1}^N \|\nabla {\widehat Y}_k^h \|_1^2 \tau\Big)^p \\
& \le C \ee \|X_0\|_1^{2p} +C \sum_{k=1}^N (1+\ee \|X(t_{k-1})\|_1^{2p(q+1)}) \tau
+ C \sum_{k=1}^N \ee \|{\widehat Y}_{k-1}^h\|_1^{2p} \tau.
\end{align*}
By the estimate \eqref{reg-u} and a discrete Gronwall inequality, we conclude \eqref{reg-aux}.
\qed   
\end{proof}

\begin{remark}
If we consider \eqref{f-coe} rather than the coercivity condition in \eqref{f-mon}, one can get an estimate of the term $\ee (\sum_{1 \le n \le N} \|\nabla Y_n\|^2 \tau)^p$ for any $p \ge 2$:
\begin{align*}
& \ee \sup_{1 \le n \le N} \|Y_n\|^{2p}
+ \ee \Big(\sum_{1 \le n \le N} \|\nabla Y_n\|^2 \tau \Big)^p
+ \ee \Big(\sum_{1 \le n \le N} \|Y_n\|_{L_\xi^{q+2}}^{q+2} \tau \Big)^p \\
& \le e^{C N \tau} (1+ \ee \|X_0\|^{2p}), \quad n \in \nn_+.
 	\end{align*}
We note that such an estimate was used in \cite[Theorem 4.1]{HS23} to control the $L_\xi^6$-norm of a tamed scheme applied to 1D SACE in deriving the overall strong convergence rate. 
Our method directly establishes the $\dot {H} ^1$-norm of the auxiliary process associated with the proposed drift-tamed-FEM, eliminating the dimension restriction.    
\end{remark}

\subsection{Strong Error Estimate}

In this part, we establish the convergence of the drift-tamed-FEM\eqref{tfem-} with a finite-time strong error estimate towards Eq. \eqref{see}. 

Denote by $E_{h,\tau}(t)=S(t)-S_{h,\tau}^n  \PP_h$ for $t\in (t_{n-1}, t_n ]$ with $1 \le n \le N$.
We will use the following estimates of $E_{h,\tau}$ in the error analysis of the tamed-FEM \eqref{tfem-}; see, e.g., \cite[Lemma 4.1]{LQ21}.

\begin{lemma} \label{lm-eht}
Let $t>0$. There exists a constant $C$ such that  
\begin{align}
& \|E_{h,\tau} (t) x\|
\le C (h^\mu+\tau^\frac\mu2) t^{-\frac{\mu-\nu}2} \|x\|_\nu,
\quad 0\le \nu\le \mu\le 2, ~ x\in \dot H^\nu, \label{eht1} \\
& \|E_{h,\tau}(\cdot) x\|_{L_t^2 L_\xi^2} 
\le C (h^{1+\mu}+\tau^\frac{1+\mu}2) \|x\|_\mu,
\quad 0\le \mu\le 1, ~ x\in \dot H^\mu. \label{eht2}
\end{align} 
\end{lemma}

Next, we establish the strong error estimate between the exact solution $X$ of Eq.~\eqref{see} and the auxiliary process $\{\widehat Y_n^h\}$ defined in \eqref{aux}. For simplicity, we assume that the initial datum $X_0$ is a deterministic element of $\dot H^1$ for the remainder of the paper. Nevertheless, all results in the next two sections remain valid for $\mathcal{F}_0$-measurable $X_0$ with sufficiently bounded higher-order moments.

\begin{lemma} \label{lm-aux} 
Let $X_0 \in \dot H^{1+\gamma}$ with $\gamma\in [0,1]$ and Assumptions \ref{ap-f}-\ref{ap-g}(1) (and (2) if $\gamma=1$) hold. 
For any $p\ge 1$, there exists a constant $C$ such that   
\begin{align} \label{err-aux}
& \sup_{1 \le n \le N} \|X(t_n)-{\widehat Y}_n^h\|_{L_\omega^p L_\xi^2}
 \le C e^{CT} (1 + \|X_0\|_{1+\gamma}^{2q+1}) (h^{1+\gamma}+\tau^{1/2} ).
\end{align}   
\end{lemma}

\begin{proof}
Let $p \ge 2$ and $n \in \nn_+$.
Subtracting the auxiliary process $\widehat Y_n^h$ defined in \eqref{aux} from the mild formulation \eqref{mild} at $t = t_n$, we obtain
\begin{align} \label{j}
J^n  :=\|X(t_n )-{\widehat Y}_n^h \|_{L_\omega^p L_\xi^2}
\le \sum_{i=1}^3 J^n_i,
\end{align}
where 
\begin{align*} 
J^n_1 & =\| E_{h,\tau}(t_n ) X_0\|_{L_\omega^p L_\xi^2}, \\
J^n_2 & =\|\sum_{i=0}^{n-1} \int_{t_i}^{t_{i+1}} [S(t_n-r) F(X)-S_{h,\tau}^{n-i} \PP_h F_\tau(X(t_i))] {\rm d}r \|_{L_\omega^p L_\xi^2},   \\
J^n_3 & =\|\sum_{i=0}^{n-1} \int_{t_i}^{t_{i+1}} [S(t_n-r) G(X)-S_{h,\tau}^{n-i} \PP_h G(X(t_i))] {\rm d}W(r) \|_{L_\omega^p L_\xi^2}.
\end{align*}
In the sequel, we treat the above three terms one by one.

The estimate \eqref{eht1} with $\mu=\nu=1+\gamma$ yields that 
\begin{align}  \label{j1}
J^n _1
\le C (h^{1+\gamma}+\tau^\frac{1+\gamma}2 ) \|X_0\|_{1+\gamma},
\quad \gamma\in [0,1].
\end{align}
To deal with the second term, we decompose it into the following three terms by Minkowski inequality:
\begin{align*} 
J^n_2 
&\le \sum_{i=0}^{n-1} \int_{t_i}^{t_{i+1}} 
\|S(t_n-r) [F(X)-F(X(t_i))] \|_{L_\omega^p L_\xi^2} {\rm d}r  \\
&\quad + \sum_{i=0}^{n-1} \int_{t_i}^{t_{i+1}} 
\|E_{h,\tau}(t_n-r) F(X(t_i))] \|_{L_\omega^p L_\xi^2} {\rm d}r \\
&\quad + \sum_{i=0}^{n-1} \int_{t_i}^{t_{i+1}} 
\|S_{h,\tau}^{n-i} [F(X(t_i))-F_\tau(X(t_i))]\|_{L_\omega^p L_\xi^2} {\rm d}r 
=: \sum_{i=1}^3 J^n_{2i}.
\end{align*} 
By \eqref{ana} with $(\mu, \nu)=(1/2, 0)$ and  the dual estimate \eqref{F-}, we have  
\begin{align*} 
J^n_{21} 
&\le C \sum_{i=0}^{n-1} \int_{t_i}^{t_{i+1}} (t_n-r)^{-1/2}
\|F(X)-F(X(t_i)) \|_{L_\omega^p \dot H^{-1}} {\rm d}r \\ 
&\le C \sum_{i=0}^{n-1} \int_{t_i}^{t_{i+1}} (t_n-r)^{-1/2}
\| (1+\|X\|^q_1+\|X(t_i)\|^q_1 ) \|X-X(t_i)\|\|_{L_\omega^p} {\rm d}r \\
& \le C \tau^{1/2} (1+\|X\|^q_{L_t^\infty L_\omega^{2pq} \dot H^1} )
\|X\|_{\CC_t^{1/2} L_\omega^{2p} L_\xi^2} \Big(\int_0^{t_n} r^{-1/2}  {\rm d}r \Big) \\  
& \le C T^{1/2} \tau^{1/2} (1+\|X\|^q_{L_t^\infty L_\omega^{2pq} \dot H^1} )
\|X\|_{\CC_t^{1/2} L_\omega^{2p} L_\xi^2}.
\end{align*}  
Using \eqref{eht1} with $(\mu, \nu)=(1+\gamma, 0)$ and the embedding \eqref{emb}, we derive  
\begin{align} \label{j22}
J^n _{22}
&\le C (h^{1+\gamma}+\tau^\frac{1+\gamma}2) 
\|F(X)\|_{L_t^\infty L_\omega^p L_\xi^2} 
\Big(\int_0^{t_n} \sigma^{-\frac{1+\gamma}2} {\rm d}\sigma\Big)  \nonumber \\
&\le C T^{\frac{1-\gamma}2} (h^{1+\gamma}+\tau^\frac{1+\gamma}2) 
(1+\|X\|_{L_t^\infty L_\omega^p \dot H^1}^{q+1}),
\quad \gamma\in [0,1).
\end{align}
The stability property \eqref{sht1} and the estimates \eqref{F-Ftau} and \eqref{reg-u} imply 
\begin{align} \label{j23}  
J^n_{23} 
& \le \sum_{i=0}^{n-1} \int_{t_i}^{t_{i+1}} 
\|F(X(t_i))-F_\tau(X(t_i))\|_{L_\omega^p L_\xi^2} {\rm d}r \nonumber  \\
&\le CT \tau^{1/2} (1+\|X\|_{L_t^\infty L_\omega^{p(2q+1)} \dot H^1}^{2q+1}).
\end{align} 

For $\gamma\in [0,1)$, combining the above three estimates and Lemma \ref{lm-u} implies 
\begin{align}  \label{j2}
J^n _2
& \le C (1+T) (h^{1+\gamma}+\tau^{1/2}) 
[(1+\|X\|^q_{L_t^\infty L_\omega^{2pq} \dot H^1} )
\|X\|_{\CC_t^{1/2} L_\omega^{2p} L_\xi^2} \nonumber \\
& \quad + 1+\|X\|_{L_t^\infty L_\omega^p \dot H^1}^{q+1}
+ \|X\|_{L_t^\infty L_\omega^{p(2q+1)} \dot H^1}^{2q+1}]  \nonumber \\ 
& \le C e^{CT} (h^{1+\gamma}+\tau^{1/2} )  
(1+ \|X_0\|_1^{2q+1}).
\end{align}

The last term $J^n _3$ can be controlled by using the continuous and discrete BDG inequalities (see, e.g., \cite[(2.22)-(2.23)]{LQ21}) as  
\begin{align*} 
(J^n _3)^2  
& \le C \sum_{i=0}^{n-1} \int_{t_i}^{t_{i+1}} 
\|S(t_n-r) G(X)-S_{h,\tau}^{n-i} \PP_h G(X(t_i))] \|_{L_\omega^p \LL_2^0}^2 {\rm d}r  \\
& \le C \sum_{i=0}^{n-1} \int_{t_i}^{t_{i+1}} 
\|S_{h,\tau}^{n-i} \PP_h [G(X)-G(X(t_i))]\|_{L_\omega^p \LL_2^0}^2 {\rm d}r  \\
& \quad + C \sum_{i=0}^{n-1} \int_{t_i}^{t_{i+1}} 
\|E_{h,\tau}(t_n-r) G(X)] \|_{L_\omega^p \LL_2^0}^2 {\rm d}r
:=\sum_{i=1}^2 (J^n _{3i})^2.
\end{align*}
Using the stability property \eqref{sht1}, conditions \eqref{g-lip}-\eqref{g1}, the estimate \eqref{eht1} with $(\mu,\nu)=(1+\gamma,1)$, and estimates \eqref{reg-u}-\eqref{hol-u}, we get 
\begin{align}  \label{j3}  
J^n _3
&\le C T^{1/2} \tau^{1/2} \|X\|_{\CC_t^{1/2} L_\omega^{2p} L_\xi^2}
+C T^{\frac{1-\gamma}2} (h^{1+\gamma}+\tau^\frac{1+\gamma}2) 
(1+\|X\|_{L_t^\infty L_\omega^p \dot H^1}) \nonumber \\ 
& \le C e^{CT} (1+\|X_0\|_1^{q+1} ) (h^{1+\gamma}+\tau^{1/2}), 
\quad \gamma\in [0,1).
\end{align} 
Putting the estimates \eqref{j1}, \eqref{j2}, and \eqref{j3} together results in
\begin{align*} 
J^n 
\le C e^{CT} (1+ \|X_0\|_1^{2q+1}) (h^{1+\gamma}+\tau^{1/2} ), 
\quad \gamma\in [0,1).
\end{align*}
This inequality, in combination with the estimate on $J^0$ that  
\begin{align*} 
J^0=\|X_0-X_0^h\|_{L_\omega^p L_\xi^2}
\le C h^{1+\gamma} \|X_0\|_{1+\gamma},
\quad \gamma\in [0,1],
\end{align*}
completes the proof of \eqref{err-aux} with $\gamma\in [0,1)$. 

To show \eqref{err-aux} for $\gamma=1$, we just need to give refined estimates for $J^n _{22}$ and $J^n _{32}$ provided $X_0 \in \dot H^2$ and \eqref{g2} holds.
Applying Minkowski inequality and using \eqref{eht1} with $\mu=2$ and $\nu=\beta\in (0,1/2)$, the embedding $\dot H^2 \hookrightarrow \dot H^\beta \cap L_\xi^\infty$, the estimate \eqref{reg-u}, and the fact  (see \cite[Lemma 3.2]{LQ21}) 
\begin{align*}
\|F(u)\|_\beta
&\le C (1+\|u\|^{q+1}_{L_\xi^\infty} +\|u\|^{q+1}_\beta ),
\quad \forall ~ u\in \dot H^\beta\cap L_\xi^\infty,
\end{align*}
we derive  
\begin{align} \label{j22+}
J^n _{22}
&\le C T^{\beta/2}(h^2+\tau) 
(1+\|X\|^{q+1}_{L_t^\infty L_\omega^p \dot H^\beta}
+\|X\|^{q+1}_{L_t^\infty L_\omega^p L_\xi^\infty})  \nonumber \\
&\le C e^{CT}(1+\|X_0\|_2^{(q+1)^2}) (h^2+\tau).
\end{align}
On the other hand, \eqref{eht1} with $(\mu, \nu)=(2, 1+\theta)$, the condition \eqref{g2} and the estimate \eqref{reg-u} imply that 
\begin{align}  \label{j32} 
J^n _{32}
& \le C T^{\theta/2} (h^2+\tau) \|G(X)\|_{L_t^\infty L_\omega^p \LL_2^{1+\theta}}   \nonumber \\
& \le C T^{\theta/2} (h^2+\tau) (1+\|X\|_{L_t^\infty L_\omega^p \dot H^{1+\theta}})  \nonumber \\ 
& \le C e^{CT} (1+\|X_0\|^{q+1}_2 ) (h^2+\tau ).
\end{align}
This completes the proof.
\qed
\end{proof}

Combining Lemma~\ref{lm-aux} with the variational approach used in the proof of Theorem~\ref{tm-lya}, we obtain the following strong convergence rate between the solution $X$ of \eqref{see} and the drift-tamed-FEM $\{Y_n^h\}$ defined in \eqref{tfem-}.

 \begin{theorem}  \label{tm-err}
Let $X_0 \in \dot H^{1+\gamma}$ with $\gamma\in [0,1]$ and Assumptions \ref{ap-f}-\ref{ap-g}(1) (and (2) if $\gamma=1$) hold. 
For any $p\ge 1$, there exists a constant $C$ such that  
\begin{align} \label{err} 
& \sup_{1 \le n \le N} \|X(t_n) - Y_n^h \|_{L_\omega^{2p} L_\xi^2} 
\le C e^{CT} (1 +\|X_0\|_{1+\gamma})^{q^2+3q+1} (h^{1+\gamma}+\tau^{1/2} ).
\end{align} 
 \end{theorem}

 \begin{proof} 
Let $1 \le n \le N$ and denote  
$e_n^h: ={\widehat Y}_n^h -Y_n^h$. 
Then $e_n^h  \in V_h$
with vanishing initial datum $e^h_0=0$.
In terms of \eqref{aux} and \eqref{full}, we have 
\begin{align*} 
{\widehat Y}_n^h 
&={\widehat Y}_{n-1}^h+\tau A_h {\widehat Y}_n^h 
+\tau \PP_h F_\tau(X(t_{n-1}))
+\PP_h G(X(t_{n-1})) \delta_{n-1} W, \\
Y_n^h  
&=Y_{n-1}^h +\tau A_h Y_n^h  
+\tau \PP_h F_\tau(Y_{n-1}^h)
+\PP_h G(Y_{n-1}^h) \delta_{n-1} W.
\end{align*}
Consequently, 
\begin{align*}
e_n^h  -e_{n-1}^h 
&=\tau A_h e_n^h  
+\tau \PP_h (F_\tau(X(t_{n-1}))-F_\tau(Y_{n-1}^h)) \nonumber \\
&\quad +\PP_h (G(X(t_{n-1}))-G(Y_{n-1}^h)) \delta_{n-1} W.
\end{align*}
Testing both sides of the above equation with $e_n^h$ and applying the elementary equality \eqref{ab}, the estimate \eqref{Ftau-mon}, and condition \eqref{g-lip}, we obtain
\begin{align*}
& \|e_n^h  \|^2-\|e_{n-1}^h\|^2 + \|e_n^h -e_{n-1}^h\|^2
+2\tau \|\nabla e_n^h  \|^2 \nonumber \\
&= 2\tau ~ {_{-1}}\<F_\tau(X(t_{n-1}))-F_\tau({\widehat Y}_{n-1}^h), e_n^h\>_1 \\
& \quad +2\tau ~{_{-1}}\<F_\tau({\widehat Y}_{n-1}^h)-F_\tau(Y_{n-1}^h), e_n^h -e_{n-1}^h\>_1 \nonumber \\
&\quad + 2\tau ~ {_{-1}}\<F_\tau({\widehat Y}_{n-1}^h)-F_\tau(Y_{n-1}^h), e_{n-1}^h \>_1 \\
&\quad + 2\<e_n^h  -e_{n-1}^h, (G(X(t_{n-1}))-G(Y_{n-1}^h))\delta_{n-1} W\>+ R_{n-1}^h \\
& \le C\tau  \|e_{n-1}^h\|^2+ \tau \|\nabla e_n^h  \|^2
+ \tau \|F_\tau(X(t_{n-1}))-F_\tau({\widehat Y}_{n-1}^h)\|_{-1}^2 \\
& \quad + 4 \tau^2 \|F_\tau({\widehat Y}_{n-1}^h)-F_\tau(Y_{n-1}^h)\|^2 + \|e_n^h  -e_{n-1}^h \|^2 \\ 
& \quad + 4 \|(G(X(t_{n-1}))-G(Y_{n-1}^h))\delta_{n-1} W\|^2
+ R_{n-1}^h,
\end{align*}  
where $R_{n-1}^h:=2\<e_{n-1}^h , (G(X(t_{n-1}) )-G(Y_{n-1}^h))\delta_{n-1} W\>$.  
It follows from the estimates \eqref{ftau'-grow} and \eqref{Ftau-} that 
\begin{align*}  
& \|e^h_n\|^2 -(1+C \tau) \|e^h_{n-1}\|^2 + \tau \|\nabla e_n^h  \|^2 \\ 
& \le \tau \|F_\tau(X(t_{n-1}))-F_\tau({\widehat Y}_{n-1}^h)\|_{-1}^2
+ 4 \tau^2 \|F_\tau({\widehat Y}_{n-1}^h)-F_\tau(Y_{n-1}^h)\|^2 \\
& \quad  + 4 \|(G(X(t_{n-1}))-G(Y_{n-1}^h))\delta_{n-1} W\|^2
+ R_{n-1}^h \\
& \le C \tau \|e_{n-1}^h\|^2
+ \tau \|F_\tau(X(t_{n-1}))-F_\tau({\widehat Y}_{n-1}^h)\|_{-1}^2 \\
& \quad   
+ 4 \|(G(X(t_{n-1}))-G(Y_{n-1}^h))\delta_{n-1} W\|^2
+ R_{n-1}^h.
\end{align*} 

Summing up the above inequality yields that  
\begin{align*}  
& \|e^h_n\|^2 + \tau \sum_{k=1}^n \|\nabla e_k^h  \|^2 \\
& \le C \tau \sum_{0 \le k \le n-1} \|e^h_k\|^2
+ \tau \sum_{0 \le k \le n-1} \|F_\tau(X(t_k))-F_\tau({\widehat Y}_k^h)\|_{-1}^2 \\
& \quad + 4 \sum_{0 \le k \le n-1} \|(G(X(t_k))-G(Y_k^h))\delta_k W\|^2 + \sum_{0 \le k \le n-1} |R_k^h|. 
\end{align*}
As in \eqref{est-sto}, we apply discrete BDG, Cauchy--Schwarz, and Young inequalities, to obtain 
\begin{align*}  
& \ee\Big( \sum_{0 \le k \le n-1} \|(G(X(t_k))-G(Y_k^h))\delta_k W\|^2\Big)^p + \ee \sup_{1 \le n \le N} \Big|\sum_{0 \le k \le n-1} R_k^h \Big|^p \\
& \le C \tau \sum_{0 \le k \le n-1} \ee \|e^h_k\|^{2p} 
+ C \sup_{0 \le k \le n-1} \ee \|X(t_k)-{\widehat Y}_k^h\|^{2p}.
\end{align*}  

Then, we use the elementary inequality \eqref{in-sum} to derive 
\begin{align*}  
\ee \sup_{1 \le n \le N} \|e^h_n\|^{2p}  
& \le C \tau \sum_{0 \le k \le n-1} \ee \|e^h_k\|^{2p} 
+ C \sup_{0 \le k \le n-1} \ee \|X(t_k)-{\widehat Y}_k^h\|^{2p} \\
& \quad + C \sup_{0 \le k \le n-1} \ee \|F_\tau(X(t_k))-F_\tau({\widehat Y}_k^h)\|_{-1}^{2p}.
\end{align*}
By discrete Gr\"onwall inequality and the estimate \eqref{Ftau-}, we obtain
\begin{align}  \label{err-xy}
 \ee \sup_{1 \le n \le N} \|e^h_n\|^{2p}  
& \le C e^{CT} \sup_{0 \le k \le n-1} (\ee \|X(t_k)-{\widehat Y}_k^h\|^{4p})^{1/2} \nonumber \\
& \qquad \times (1+ \|X\|_{L_t^\infty L_\omega^{4p q} \dot H^1}^{2pq} 
+\sup_{0 \le n \le N-1} \|{\widehat Y}_k^h\|_{L_\omega^{4p q} \dot H^1}^{2pq}).  
\end{align} 
Then, we conclude \eqref{err} from \eqref{err-aux}, \eqref{reg-u}, \eqref{reg-aux}, and \eqref{err-xy}.
\qed
 \end{proof}

\section{High-order Convergence Rate in Additive Noise Case}
\label{sec5}

In this section, we analyze the strong convergence rates in the additive noise case, i.e., when the diffusion coefficient of Eq.~\eqref{see} is a constant operator $G(\cdot) \equiv G$: 
\begin{align} \label{see-}
{\rm d} X(t)=(A X(t)+F(X(t))) {\rm d}t+ G{\rm d}W, \quad  t \ge 0;
\quad X(0)=X_0.
\end{align}
In this case, conditions \eqref{g-lip}-\eqref{g1} are equivalent to the assumption $G \in \mathcal{L}_2^1$, which was imposed in \cite{QW19} as $\|(-A)^{1/2}\|_{\mathcal{L}_2^0} < \infty$ with $g(\cdot) \equiv 1$.
 
We consider the drift-tamed-FEM \eqref{tfem-} applied to Eq.~\eqref{see} with a constant diffusion operator $G(\cdot) \equiv G$ and with initial condition $Y_0^h = \mathcal{P}_h X_0$ and 
 \begin{align} \label{tem*}
Y_n^h=Y_{n-1}^h + A_h Y_n^h \tau
+ \PP_h F_\tau (Y_{n-1}^h) \tau
 	+ \PP_h G \delta_{n-1} W,
\quad n \in \nn_+. 
 \end{align} 
Our main purpose is to establish a high-order convergence result under the additional polynomial growth condition on the second derivative $f''$ of $f$:
 \begin{align}\label{f''}
|f''(\xi)| \le L_3' +L_4' |\xi|^{\max\{0, q-1\}},
\quad \xi \in \rr,
  \end{align}
  for some positive constants $L_3'$ and $L_4'$.
It is clear that the condition \eqref{f''} is consistent with Assumption \ref{ap-f}. 

Now, we can state and prove the higher-order convergence for the drift-tamed-FEM \eqref{tem*} applied to Eq. \eqref{see-}.

 \begin{theorem}  \label{tm-err-} 
Let $X_0 \in \dot H^{1+\gamma} \cap L_\xi^{2(3q+1)}$ with $\gamma\in [0,1]$, and $G \in \LL_2^1$ (and there exists a positive constant $\theta$ such that $G \in \LL_2^{1+\theta}$ if $\gamma=1$).
Under the condition \eqref{f''}  and Assumption \ref{ap-f}, for any $p\ge 1$, there exists a constant $C$ such that  
\begin{align} \label{err-} 
& \sup_{1 \le n \le N} \|X(t_n) - Y_n^h \|_{L_\omega^{2p} L_\xi^2} \nonumber \\
& \le C e^{CT} (1 +\|X_0\|_{1+\gamma}^{q^2+3q+1}
+ \|X_0\|_{L_\xi^{2(3q+1)}}^{3q+1}) (h^{1+\gamma}+\tau^\frac{1+\gamma}2 ).
\end{align} 
 \end{theorem}

\begin{proof}
Let $p \ge 2$ and $\gamma\in [0,1]$.
As \eqref{err-xy} in the proof in Theorem \ref{tm-err}, we have 
\begin{align} \label{err-xy+}  
& \ee \sup_{1 \le n \le N} \|X(t_n) - Y_n^h\|_{L_\omega^{2p} L_\xi^2}   \nonumber \\ 
& \le C e^{CT} \sup_{0 \le k \le n-1} \|X(t_k)-{\widehat Y}_k^h\|_{L_\omega^{4p} L_\xi^2} \nonumber \\
& \qquad \times (1+ \|X\|_{L_t^\infty L_\omega^{4p q} \dot H^1}^q 
+\sup_{0 \le n \le N-1} \|{\widehat Y}_k^h\|_{L_\omega^{4p q} \dot H^1}^q). 
\end{align} 
Thus, it suffices to establish the following refined estimate between $X$ and $\widehat Y_n^h$ when $G(\cdot) \equiv G$. Following the proof of Lemma~\ref{lm-aux}, we only need to provide refined estimates for $J_{21}^n$, $J_{23}^n$, and $J_3^n$ under the additional condition \eqref{f''}.

To show a refined estimate of $J^n_{21}$, note that for $r\in [t_i,t_{i+1})$,
\begin{align*} 
X(t_{i+1})
&=S(t_{i+1}-r) X(r)
+\int_r^{t_{i+1}} S(t_{i+1}-\sigma) F(X(\sigma)) {\rm d}\sigma \\
& \quad +\int_r^{t_{i+1}} S(t_{i+1}-\sigma) G {\rm d}W(\sigma).
\end{align*}
Then, using Taylor formula leads to the splitting of $J^n_{21}$ into 
\begin{align*} 
J^n_{21} 
&=\|\sum_{i=0}^{n-1} \int_{t_i}^{t_{i+1}} S(t_n-r) 
[F(X(r))-F(X(t_{i+1}))] {\rm d}r \|_{L_\omega^p L_\xi^2} \\
&\le \Big \|\sum_{i=0}^{n-1} \int_{t_i}^{t_{i+1}} S(t_n-r) 
DF(X(r)) (S(t_{i+1}-r)-{\rm Id}) X(r) {\rm d}r \Big\|_{L_\omega^p L_\xi^2} \\
&\quad + \Big \|\sum_{i=0}^{n-1} \int_{t_i}^{t_{i+1}} S(t_n-r) 
DF(X(r)) \Big(\int_r^{t_{i+1}} S(t_{i+1}-\sigma) F(X(\sigma)) {\rm d}\sigma\Big) {\rm d}r \Big\|_{L_\omega^p L_\xi^2} \\
&\quad + \Big \|\sum_{i=0}^{n-1} \int_{t_i}^{t_{i+1}} S(t_n-r) 
DF(X(r)) \Big(\int_r^{t_{i+1}} S(t_{i+1}-\sigma) G {\rm d}W(\sigma) \Big)  {\rm d}r \Big\|_{L_\omega^p L_\xi^2} \\
&\quad +\Big \|\sum_{i=0}^{n-1} \int_{t_i}^{t_{i+1}} S(t_n-r) 
R_F(X(r), X(t_{i+1})) {\rm d}r \Big\|_{L_\omega^p L_\xi^2}
=:\sum_{i=1}^4 J^n_{21i}, 
\end{align*}
where $R_F$ denotes the remainder term
\begin{align*} 
& R_F(X(r), X(t_{i+1})) \\
&:=\int_0^1 D^2 F (X(r)+\lambda (X(t_i)-X(r)) )
 (X(t_i)-X(r), X(t_i)-X(r)) (1-\lambda) {\rm d}\lambda \\
&=\int_0^1 f'' (X(r)+\lambda (X(t_i)-X(r)) )
 (X(t_i)-X(r))^2 (1-\lambda) {\rm d}\lambda.
\end{align*}
We shall estimate $J^n_{21i}$, $i=1,2,3,4$, successively.

Since for $\delta\in (3/2, 2)$, $\dot H^\delta\hookrightarrow \CC$, it follows by the dual argument that
\begin{align} \label{l1} 
\|x\|_{-\delta}\le C \|x\|_{L_\xi^1},
\quad x\in L_\xi^1.
\end{align}
This inequality, in conjunction with Minkowski and Cauchy--Schwarz inequalities, the condition \eqref{f'}, the embedding \eqref{emb}, and the estimates \eqref{ana} with $(\mu, \nu)=(\delta, 0)$ and \eqref{ana-hol} with $\rho=1+\gamma$, yields that 
\begin{align} \label{j211} 
J^n_{211} 
&\le C \sum_{i=0}^{n-1} \int_{t_i}^{t_{i+1}} (t_n-r)^{-\frac\delta2}
\|DF(X)\|_{L_\omega^{2p} L_\xi^2} 
\|(S(t_{i+1}-r)-{\rm Id}) X\|_{L_\omega^{2p} L_\xi^2}  {\rm d}r \nonumber \\
&\le C \tau^\frac{1+\gamma}2 T^{1-\frac \delta 2} (1+\|X\|_{L_t^\infty L_\omega^{2pq} \dot H^1}^q) \|X\|_{L_t^\infty L_\omega^{2pq} \dot H^{1+\gamma}},
\quad \gamma \in [0, 1].
\end{align}
A similar argument implies that 
\begin{align} \label{j212} 
J^n_{212} 
&\le C \sum_{i=0}^{n-1} \int_{t_i}^{t_{i+1}} \int_r^{t_{i+1}} 
(t_n-r)^{-\frac\delta2}
\|DF(X(r))\|_{L_\omega^{2p} L_\xi^2} 
\|f(X(\sigma))\|_{L_\omega^{2p} L_\xi^2} {\rm d}\sigma {\rm d}r \nonumber \\
&\le C \tau T^{1-\delta/2} (1+\|X\|_{L_t^\infty L_\omega^{2pq} \dot H^1}^{q+1}).
\end{align}

To estimate the third term $J^n_{213}$, we apply the stochastic Fubini theorem, discrete and continuous BDG, and Cauchy--Schwarz inequalities to derive 
\begin{align*}
& (J^n_{213})^2 \\
&=\Big\| \sum_{i=0}^{n-1} \int_{t_i}^{t_{i+1}} \int_{t_i}^{t_{i+1}} \chi_{[r,t_{i+1})}(\sigma) S(t_n-r) DF(X(r)) S(t_{i+1}-\sigma) G {\rm d}r {\rm d}W(\sigma) \Big\|_{L_\omega^p L_\xi^2}^2 \\ 
&\le C \tau \sum_{i=0}^{n-1} \int_{t_i}^{t_{i+1}} 
\int_{t_i}^{t_{i+1}} \|S(t_n-r) DF(X(r)) S(t_{i+1}-\sigma) G\|_{L_\omega^p \LL_2^0}^2 {\rm d}r  
{\rm d}\sigma   \nonumber \\
&\le C \tau \sum_{i=0}^{n-1} 
\Big(\int_{t_i}^{t_{i+1}} \|DF(X)\|^2_{L_\omega^p \dot H^{-1}}  {\rm d}r \Big)
\Big(\int_{t_i}^{t_{i+1}}  \|S(\sigma) G\|_{\LL_2^1}^2  {\rm d}\sigma \Big).
\end{align*}
This, in combination with \eqref{f'} and the assumption $G\in \LL_2^1$, shows that 
\begin{align} \label{j213}
J^n_{213} \le C \tau (1+\|X\|_{L_t^\infty L_\omega^{pq} \dot H^1}^q).
\end{align}
Finally, in terms of the estimates \eqref{ana} with $(\mu, \nu)=(\delta, 0)$ and \eqref{l1}, Minkowski and H\"older inequalities, the condition \eqref{f''}, and the embeddings \eqref{emb} and $\dot H^{1/2} \hookrightarrow L_\xi^3$, we estimate the last term $J^n_{214}$ similarly to $J^n_{211}$ by
\begin{align*}
J^n_{214} 
&\le C \sum_{i=0}^{n-1} \int_{t_i}^{t_{i+1}} (t_n-r)^{-\frac\delta2}
\|R_F(X, X(t_{i+1}))\|_{L_\omega^p L_\xi^1} {\rm d}r \nonumber \\
&\le C T^{1-\delta/2} \tau^{(1/2+\gamma) \wedge 1} (1+\|X\|^{q-1}_{L_t^\infty L_\omega^{3p(q-1)} L_\xi^{3(q-1)}} )\|X\|^2_{\CC_t^{[(1/2+\gamma) \wedge 1]/2} L_\omega^{3p} \dot H^{1/2}}.
\end{align*} 
It is clear that $\dot H^1 \subset L_\xi^{3(q-1)}$ as $3(q-1)\le 2(q+1)$ (since $q \le 2$) for $d=3$ and $\dot H^1 \subset L_\xi^p$ for any $p \in [1, \infty)$ for $d=1, 2$.
Consequently,   
\begin{align} \label{j214}
J^n_{214} 
\le C \tau^{(\frac12+\gamma) \wedge 1} T^{1-\frac\delta2} (1+\|X\|^{q-1}_{L_t^\infty L_\omega^{3p(q-1)} \dot H^1} )\|X\|^2_{\CC_t^{[(\frac12+\gamma) \wedge 1]/2} L_\omega^{3p} \dot H^{\frac12}}.
\end{align} 

Collecting the above four estimates \eqref{j211}-\eqref{j214} together, 
we have 
\begin{align} \label{j21}
J^n_{21} 
& \le C(1+T)  \tau^\frac{1+\gamma}2 [1+\|X\|_{L_t^\infty L_\omega^{2pq} \dot H^1}^{q+1}
+ (1+\|X\|_{L_t^\infty L_\omega^{2pq} \dot H^1}^q) \|X\|_{L_t^\infty L_\omega^{2pq} \dot H^{1+\gamma}} \nonumber  \\
& \qquad \qquad \qquad \qquad + (1+\|X\|^{q-1}_{L_t^\infty L_\omega^{3p(q-1)} \dot H^1} )\|X\|^2_{\CC_t^{[(1/2+\gamma) \wedge 1]/2} L_\omega^{3p} \dot H^{1/2}}] \nonumber \\
& \le 
\begin{cases}
C e^{CT} \tau^\frac{1+\gamma}2 (1+\|X_0\|^{3q+1}_{1+\gamma}), \quad & \gamma\in [0,1), \\
C e^{CT} \tau^\frac{1+\gamma}2 (1+\|X_0\|^{q^2+3q}_{1+\gamma}+\|X_0\|^{q^2+3q}_2), \quad & \gamma=1,
\end{cases}  
\end{align}
where in the last step we used the estimates \eqref{reg-u}-\eqref{hol-u}.

For the term $J^n_{23}$, we use the estimate \eqref{F-Ftau+}, instead of \eqref{F-Ftau}, and \eqref{reg-u+} to derive 
\begin{align} \label{j23-} 
J^n_{23} 
& \le \sum_{i=0}^{n-1} \int_{t_i}^{t_{i+1}} 
\|F(X(t_i))-F_\tau(X(t_i))\|_{L_\omega^p L_\xi^2} {\rm d}r \nonumber  \\
&\le CT \tau (1+\|X\|_{L_t^\infty L_\omega^{p(3q+1)} L_\xi^{2(3q+1)}}^{3q+1}).
\end{align}   
The last term $J^n_3$ ($J^n_{31}$ vanishes in the additive noise case) can be handled by BDG inequality and \eqref{eht2} with $\mu=1$:
\begin{align}  \label{j3+} 
J^n_3 
&=\Big\|\int_0^{t_n} 
E_{h,\tau}(t_n-r) G {\rm d}W(r) \Big \|_{L_\omega^p L_\xi^2}  \nonumber  \\
& \le C \Big(\int_0^{t_n}  \|E_{h,\tau}(\sigma) G \|^2_{\LL_2^0} {\rm d}\sigma \Big)^{1/2} \nonumber  \\
& \le C \|G \|_{\LL_2^1} (h^2+\tau).
\end{align}
Now we can conclude \eqref{err-}, taking into account \eqref{err-xy+}, the error estimates \eqref{j1}, \eqref{j21}, \eqref{j22} and \eqref{j22+},  \eqref{j23-}, and \eqref{j3+} for $J^n_1$, $J^n_{21}$, $J^n_{22}$, $J^n_{23}$, and $J^n_3$, respectively, and the moment estimates \eqref{reg-u}, \eqref{hol-u}, \eqref{reg-u+}, and \eqref{reg-aux}.
\qed
\end{proof}

 \begin{remark} \label{rk-reg}
For $d = 1$ or $d = 2$, since $\dot H^1 \subset L_\xi^p$ for any $p \in [1,\infty)$, the assumption $X_0 \in L_\xi^{2(3q+1)}$ is automatically satisfied. Consequently, the strong error estimates for the additive noise case in Theorem \ref{tm-err-} hold provided $X_0 \in \dot H^{1+\gamma}$ with $\gamma \in [0,1]$. 
The need for the $L_\xi^{2(3q+1)}$-regularity of the exact solution arises from the term $J_{23}^n$ in \eqref{j23-}; this term cannot be controlled using the $\dot H^1$-norm of the solution for the 3D SACE corresponding to  $q = 2$. 
 \end{remark}

 \section*{Declarations}

{\bf Conflict of interest}. The authors have no competing interests to declare relevant to this article's content.
The data are available from the corresponding author on reasonable request.

 \section*{Acknowledgements}

The first author is supported by the National Natural Science Foundation of China (NNSFC), No. 12671474, Guangdong Basic and Applied Basic Research Foundation, No. 2024A1515012348, and Shenzhen Basic Research Special Project (Natural Science Foundation) Basic Research (General Project), No. JCYJ20240813094919026; NNSFC, Nos. 12101296, 12371409, and W2431008 support the second author.

 \bibliographystyle{plain}
  \bibliography{bib.bib}

\end{document}